\documentclass[
  a4paper,
  11pt,
  oneside
]{article}

\usepackage[utf8]{inputenc}
\usepackage{lmodern}
\usepackage[german,english]{babel}

\usepackage[numbers]{natbib}
\usepackage{xcolor}

\usepackage{enumerate}

\usepackage{graphicx}   
\usepackage{float}      
\usepackage{wrapfig}    
\usepackage{caption}    
\usepackage{tikz}       
  \usetikzlibrary{arrows.meta}
  \usetikzlibrary{calc}

\usepackage{pgf}
\usepackage{pgfplots}
\usepackage{multirow}
\usepackage{rotating}
\usepackage{booktabs}  

\usepackage{url}
\usepackage{authblk}
\usepackage{hyperref}
  \hypersetup{
    colorlinks,
    citecolor=black,
    filecolor=black,
    linkcolor=black,
    urlcolor=black
  }

\usepackage{mathtools}  
\usepackage{xfrac}      
\usepackage{amssymb}    
\usepackage{amsmath}    
\usepackage{amsfonts}
\usepackage{units}
\numberwithin{equation}{section}
\numberwithin{figure}{section}
\numberwithin{table}{section}
\usepackage{amscd}

\usepackage{cleveref}

\usepackage{microtype}  
\usepackage{relsize}    
\usepackage[            
  top=2.5cm, 
  bottom=2.0cm, 
  left=2.5cm, 
  right=2.5cm
]{geometry}

\newtheorem{problem}{Problem}
\newtheorem{remark}{Remark}

\newcommand{\bs}[1]{{#1}}
\newcommand{\hOmega}{\widehat\Omega}

\newcommand{\hGamma}{\widehat\Gamma}

\newcommand{\hn}{\bs{\hat n}}
\newcommand{\hcalA}{\hat{\mathcal{A}}}

\newcommand{\hrho}{\bs{\hat\rho}}

\newcommand{\hx}{\bs{\hat x}}
\newcommand{\hF}{\bs{\widehat{F}}}
\newcommand{\hJ}{\hat{J}}
\newcommand{\hE}{\bs{\widehat{E}}}

\newcommand{\hnabla}{\bs{\widehat\nabla}}
\newcommand{\hSigma}{\bs{\widehat\Sigma}}

\newcommand{\hsigma}{\bs{\widehat\sigma}}

\newcommand{\hX}{\bs{\widehat{X}}}

\newcommand{\hv}{\bs{\hat v}}
\newcommand{\hu}{\bs{\hat u}}
\newcommand{\hp}{\bs{\hat p}}
\newcommand{\hw}{\bs{\hat w}}

\newcommand{\hU}{\bs{\hat U}}

\newcommand{\hPsi}{\bs{\hat\Psi}}

\newcommand{\hpsi}{\bs{\hat\psi}}

\newcommand{\belos}{\texttt{Belos}}
\newcommand{\epetra}{\texttt{Epetra}}
\newcommand{\dealii}{\texttt{deal.II}}
\newcommand{\frosch}{\texttt{FROSch}}
\newcommand{\kokkos}{\texttt{Kokkos}}
\newcommand{\kokkoskernels}{\texttt{KokkosKernels}}
\newcommand{\stratimikos}{\texttt{Stratimikos}}
\newcommand{\tpetra}{\texttt{Tpetra}}
\newcommand{\trilinos}{\texttt{Trilinos}}
\newcommand{\xpetra}{\texttt{Xpetra}}

\begin{document}

\title{%
  Coupling deal.II and FROSch: A Sustainable and Accessible (O)RAS  Preconditioner
}

\author[1]{A. Heinlein}
\author[2,3]{S. Kinnewig}
\author[2,3]{T. Wick}

\affil[1]{%
  TU Delft,
  DIAM, Faculty of EEMCS Numerical Analysis,
  Mekelweg 4,
  2628 Delft,
  Netherlands
}

\affil[2]{%
  Leibniz University Hannover,
  Institute of Applied Mathematics,
  Welfengarten 1,
  30167 Hannover,
  Germany
}

\affil[3]{%
  Cluster of Excellence PhoenixD (Photonics, Optics, and
  Engineering - Innovation Across Disciplines),
  Leibniz University Hannover,
  Germany
}

\date{}

\maketitle

\begin{abstract}
\noindent%
In this work, restricted additive Schwarz (RAS) and 
optimized restricted additive Schwarz (ORAS) preconditioners from the {\trilinos} package {\frosch} (Fast and Robust Overlapping Schwarz) are employed to solve model problems implemented using {\dealii} (differential equations analysis library). Therefore, a {\tpetra}-based interface for coupling {\dealii} and {\frosch} is implemented. While RAS preconditioners have been available before, ORAS preconditioners have been newly added to {\frosch}. The {\dealii}--{\frosch} interface works for both Lagrange-based and N\'ed\'elec finite elements. 
Here, as model problems, nonstationary, nonlinear, variational-monolithic fluid-structure interaction and the indefinite time-harmonic Maxwell's equations are considered. Several numerical experiments in two and three spatial dimensions confirm the performance of the preconditioners as well as the {\frosch}-{\dealii} interface. In conclusion, the overall software interface is 
straightforward and easy to use while giving satisfactory solver performances for challenging PDE systems.
\end{abstract}

\section{Introduction}
Domain decomposition methods (DDMs) \cite{smith_domain_1996,quarteroni_domain_1999,toselli_domain_2005,dolean_introduction_2015} are a class of solvers, and preconditioners
that allow for efficiently solving complex model problems
arising from scientific and engineering applications. Their robustness enables the application to challenging problems that would otherwise require the use of direct solvers, which is infeasible for large problems due to the superlinear complexity of direct solvers.
The restricted additive Schwarz (RAS) preconditioner~\cite{cai_restricted_1999} and optimized restricted additive Schwarz (ORAS) preconditioner~\cite{Art:Gander:06} are popular choices of one-level Schwarz domain decomposition preconditioners. While RAS preconditioners are closely related to Lions' algorithm~\cite{lions_schwarz_1988}, as discussed in~\cite{efstathiou_why_2003}, ORAS preconditioners are a typical choice of preconditioners for wave problems due to a different choice of the transmission condition on the subdomain boundaries. However, both RAS and ORAS preconditioners are inherently unsymmetric, even if the original problem is symmetric, and therefore, cannot be employed with the conjugate gradient (CG) method; here, we consider the GMRES (generalized minimal residual) method~\cite{saad_gmres_1986} instead.

One-level DDMs, such as one-level RAS and ORAS preconditioners, are generally not scalable to large numbers of subdomains; the number of iterations will increase with the number of subdomains.
To overcome this problem, a second level, which corresponds to solving a global coarse problem of a small dimension, can be introduced. If well chosen, the coarse problem enables the global transport of information required for numerical scalability, that is, scalability for increasing numbers of subdomains; however, since the coarse level introduces a global problem, it can also become a bottleneck for parallel scalability if its dimension is too large. If the coarse problem can be solved efficiently (in parallel), two-level Schwarz preconditioners are a compelling choice for high-performance computing (HPC) applications.
In order to avoid implementing RAS or ORAS preconditioners from scratch for various applications,
it is desirable to have a sustainable, accessible implementation that can be adjusted to other application problems with minimal effort. To address this, we choose the {\frosch}~\cite{heinlein_frosch_2020,heinlein_parallel_2016-1} domain decomposition solver framework, which is part of {\trilinos}, as the basis for implementing the preconditioners. Then, we 
couple it with the finite element library {\dealii}~\cite{deal2020,Software:dealii:95}. {\dealii} is notable for its active community, comprehensive documentation, and numerous online tutorial examples. 

The key contributions of this work are the implementation of the ORAS preconditioner in {\frosch} as well as the coupling of
{\frosch} and {\dealii}. To this end, we construct and implement a software
interface based on the newly developed {\tpetra}-Interface of {\dealii}. The performance is then demonstrated for
two challenging model problems. First, nonstationary, nonlinear
Fluid-Structure Interaction (FSI) is considered in both two and three spatial dimensions. A two-level RAS preconditioner with a specific coarse space from {\frosch} is employed here. Second, the new ORAS preconditioner, as well as the new {\frosch}-{\dealii} interface, are utilized to solve
the time-harmonic Maxwell's equations,
again, both in two and three spatial dimensions.
The main effort has been in the easiness of using the proposed software interface while giving satisfactory solver performances 
for challenging PDE problems, namely FSI and time-harmonic Maxwell's equations. 
The code is available on GitHub\footnote{\url{https://github.com/kinnewig/dealii-FROSch-interface}}, where we provide one example of the RAS preconditioner and one example of the ORAS preconditioner.

The outline of this paper is as follows. In~\Cref{sec:schwarz}, we provide a detailed introduction to Schwarz preconditioners,
as well as a literature overview about Schwarz preconditioners. In~\Cref{sec:dealii_and_frosch}, we describe the employed
software packages, i.e., {\frosch}, {\dealii}, and the interface between these two. The partial differential equations considered as examples are introduced in~\Cref{sec:applications}. Finally, in~\Cref{sec:results}, several numerical tests in two and three spatial dimensions for FSI and Maxwell’s equations substantiate our newly developed software interface.

\section{Overlapping Schwarz Preconditioners} \label{sec:schwarz}
In this section, we introduce different versions of overlapping 
Schwarz preconditioners.
Employing the Finite Element Method (FEM) to discretize (initial) boundary value problems on some computational domain $\Omega$ results at some point 
(after possible linearization) in sparse linear equation systems of the form
\begin{equation} \label{eq:les}
	A x = b.
\end{equation}
For future reference, we use $h$ to represent the mesh parameter, and the finite element space is denoted as $V \coloneqq V^h(\Omega)$.

\begin{figure}[t]
  \begin{minipage}{0.45\textwidth}
	  \centering
    \includegraphics[width=0.7\textwidth]{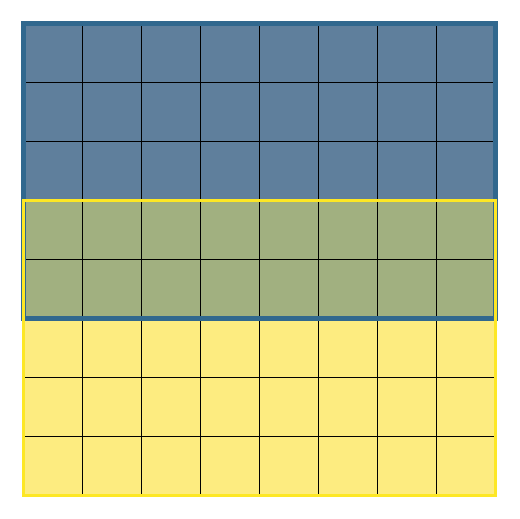}
  \end{minipage}
  \hfill
  \begin{minipage}{0.45\textwidth}
	  \centering
    \includegraphics[width=0.7\textwidth]{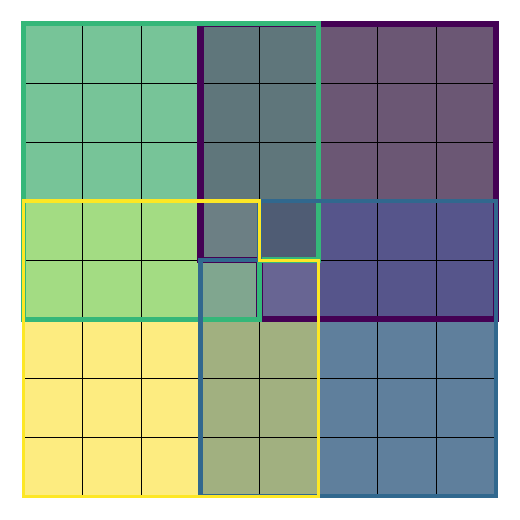}
  \end{minipage}
	\caption{
  	\label{fig:overlapping_dd}
      \textbf{Left:} A square domain divided into two subdomains, where the size of both subdomains was increased by one layer. 
      \textbf{Right:} A square domain divided into four subdomains, where the size of all subdomains was increased by one layer. 
      In both cases, the overlap was computed based on the dual graph.
    }
\end{figure}

\subsection{Additive Schwarz preconditioners} \label{sec:as}
We decompose $\Omega$ into $N$ non-overlapping subdomains $\Omega_1,\ldots,\Omega_N$. In practice, this decomposition may, for 
instance, be obtained via a geometrical approach or via partitioning the dual graph of the triangulation using a graph partitioning 
tool, such as \texttt{METIS}~\cite{karypis_fast_1998}, \texttt{Zoltan}~\cite{boman_zoltan_2012}, or \texttt{p4est}~\cite{Software:p4est:11}. In the 
dual graph, the nodes of the graph are the finite elements, and two elements are connected via an edge if the finite elements have a
common edge in the triangulation. Geometrical approaches, such as decomposing a rectangular domain into rectangular subdomains, are 
primarily limited to simple geometries and structured meshes, whereas graph partitioning approaches also apply to unstructured meshes.

Next, we extend these subdomains by $k$ layers of elements, resulting in an overlapping domain decomposition with the subdomains $\Omega_1',
\ldots,\Omega_N'$; see~\Cref{fig:overlapping_dd}. 
We denote the size of the overlap by $\delta = kh$ and define the local finite element spaces $V_i$ on the local overlapping subdomains. Due to our construction, 
it follows that $V_i \subset V$, and we can 
define restriction operators $R_i: V \to V_i$ and corresponding prolongation operators $P_i: V_i \to V$, for all $1 \leq i \leq N$. Both can be
represented as sparse binary matrices, and $P_i = R_i^\top$ is the transposed operator of $R_i$. Then, the additive overlapping Schwarz (AS)
preconditioner 
is formally written as
\begin{equation} \label{eq:as}
	M_{\mathrm AS}^{-1} 
	=
	\sum_{i=1}^N 
	P_i A_i^{-1} R_i.
\end{equation}
In the case of local exact solvers, we define the local subdomain matrices as 
\begin{equation} \label{eq:ai}
	A_i = R_i A P_i, \quad\forall 1 \leq i \leq N.
\end{equation}
Therefore, we choose $V_i = V_0^h (\Omega_i')$, that is, the local finite element spaces are the restrictions of the global finite element space $V$ to the local subdomains and
with homogeneous Dirichlet boundary conditions on $\partial\Omega_i'$. Then, $A_i$ can be obtained from $A$ by extracting the submatrix
corresponding to the interior finite element nodes of $\Omega_i'$. This also directly implies that the $A_i$ are invertible. For more
details, we refer to the standard domain decomposition 
literature~\cite{smith_domain_1996,quarteroni_domain_1999,toselli_domain_2005,dolean_introduction_2015}.

For a Laplacian model problem, the matrix $A$ is symmetric positive definite, and we can solve~\Cref{eq:les} using the preconditioned conjugate gradient (PCG) method
with the preconditioner given in~\Cref{eq:as}; we notice that $M_{\rm AS}^{-1}$ is then also symmetric positive definite. The convergence of
PCG can then be bounded in terms of $\kappa$, namely the condition number of the preconditioned system matrix $M_{\rm AS}^{-1} A$, i.e.,
\begin{equation} \label{eq:cond:as}
	\kappa (M_{\rm AS}^{-1} A) 
	\leq 
	C \left( 1 + \frac{1}{H \delta} \right).
\end{equation}
Here, $H$ is the subdomain diameter, and $C$ is a constant that is independent of geometric constants $H$, $h$, and $\delta$; see, for
instance~\cite{smith_domain_1996,quarteroni_domain_1999}.

\begin{remark} \label{remark:algebraic}
  The aforementioned approaches for the construction of overlapping domain decomposition cannot be performed in a fully algebraic
  way. We refer to \textit{fully algebraic} as a method that requires only the fully assembled system matrix $A$. In particular, the
  geometry and/or the triangulation are required. The graph partitioning-based approach can also be carried out based on the sparsity
  pattern of $A$, that is, based on the node graph of the triangulation. Then, the overlap can be obtained by adding $k$ layers of
  adjacent finite element nodes, where adjacency is defined based on the node graph. The resulting overlapping domain decomposition method
  then slightly differs, but the preconditioner becomes fully algebraic.
\end{remark}

\subsection{Restricted additive Schwarz preconditioners} \label{sec:ras}
In~\cite{cai_restricted_1999}, Cai and Sarkis introduced a restricted additive Schwarz (RAS) preconditioner, which often converges faster than the standard AS preconditioner and has a reduced communication cost in a distributed memory parallelization. However, the preconditioner is non-symmetric and, therefore, requires a different iterative solver than PCG, such as the generalized minimal residual (GMRES) method~\cite{saad_gmres_1986}. It is obtained by keeping $A_i$ and $R_i$ in~\Cref{eq:as} while replacing the prolongation operators $P_i$ by $\hat P_i$:
\begin{equation} \label{eq:ras1}
	M_{\rm RAS}^{-1} 
	=
	\sum_{i=1}^N 
	\hat P_i A_i^{-1} R_i.
\end{equation}
Here, we define $\hat P_i$ in such a way that
\begin{equation} \label{eq:pou1}
	\sum_{i=1}^N 
	\hat P_i R_i
	= 
	I.
\end{equation}
In the classical RAS preconditioner, this is achieved by partitioning the finite element nodes of the overlapping subdomains $\Omega_1',\ldots,\Omega_N'$ in a unique way. As discussed by Efstathiou and Gander in~\cite{efstathiou_why_2003}, this preconditioner relates to 
Lions' classical Schwarz method~\cite{lions_schwarz_1988}. \Cref{eq:pou1} 
can alternatively be written
using diagonal scaling matrices $D_i$:
\begin{equation} \label{eq:pou2}
	\sum_{i=1}^N 
	P_i D_i R_i
	= 
	I.
\end{equation}
The previous choice of $\hat P_i$ corresponds to using binary scaling matrices $D_i$. However, we can also choose $D_i$ to scale each finite element node by its inverse multiplicity with respect to the overlapping domain decomposition into $\Omega_1',\ldots,\Omega_N'$. Then, we obtain the following alternative formulation of the RAS preconditioner
\begin{equation} \label{eq:ras2}
	M_{\rm RAS}^{-1} 
	=
	\sum_{i=1}^N 
	P_i D_i A_i^{-1} R_i.
\end{equation}
When using a non-binary scaling $D_i$, the preconditioner is sometimes called a scaled additive Schwarz (SAS) preconditioner.

\subsection{Optimized Schwarz preconditioners} 
\label{sec:oas}
Finally, we introduce optimized additive Schwarz (OAS) and optimized restricted additive Schwarz (ORAS) 
preconditioners~\cite{st-cyr_optimized_2007}. These are often employed for wave-type problems,
such as Helmholtz equations or 
Maxwell's equations; cf.~\Cref{sec:maxwell}. Here, in particular, standard AS preconditioners often lead to suboptimal
convergence rates compared to optimized Schwarz preconditioners.
To this end,
let us consider that~\Cref{eq:les} results from the discretization of the Helmholtz equation
$$
- \Delta u - \omega^2 u = f \quad \text{in } \Omega, 
$$
with some boundary conditions on $\partial\Omega$. Here, $\omega$
is the so-called wavenumber, and $f$ is a source term.

In optimized Schwarz preconditioners, we replace the local problem matrices $A_i$ in~\Cref{eq:as} with the matrices $B_i$ corresponding 
to local problems with different boundary conditions; recall 
that we automatically obtained homogeneous Dirichlet boundary conditions by the algebraic construction in~\Cref{eq:ai}.
In particular, the OAS preconditioner reads 
\begin{equation} \label{eq:oas}
	M_{\rm OAS}^{-1} 
	=
	\sum_{i=1}^N 
	P_i B_i^{-1} R_i,
\end{equation}
and the ORAS preconditioner reads
\begin{equation} \label{eq:oras}
	M_{\rm ORAS}^{-1} 
	=
	\sum_{i=1}^N 
	P_i D_i B_i^{-1} R_i.
\end{equation}
In the most simple case, the matrix $B_i$ is the discretization of the differential operator $- \Delta u - \omega^2 u$
using the local finite element space $V_i$ and with Robin boundary operator
$\left( \frac{\partial}{\partial_n} + \alpha \right)$, $\alpha \in \mathbb{R}$,
on the nodes of $\partial\Omega_i' \setminus \partial\Omega$. The Robin parameter $\alpha$ is a hyperparameter of the method, which can be tuned and gives the name to the optimized Schwarz method. In practice, the optimization can be done either theoretically~\cite{Art:Gander:06} or computationally via numerical calibration. In a Krylov iteration, the application of $B_i^{-1}$ typically corresponds to solving a system of the form
$$
	B_i y_i = R_i (b - A x^{(k)}),
$$
for $y_i$. The right-hand side is a restriction of the residual to $\Omega_i'$, and $x^{(k)}$ is the iterate in the $k$th iteration. In particular, the restriction $R_i$ is chosen such that $R_i (b - A x^{(k)}) = 0$ for the nodes on $\partial\Omega_i' \setminus \partial\Omega$. As a result, we enforce the Robin boundary condition
\begin{equation} \label{eq:robin}
	\left( \frac{\partial}{\partial_n} + \alpha \right) = 0
\end{equation}
on those nodes.

Other boundary conditions employing a perfectly matched
layer (PML)~\cite{Art:Dolean:24::D27} or other types of absorbing boundary conditions~\cite{Art:Gander:16}, are more complex and depend on more than one parameter but may yield even better convergence than Robin boundary conditions. 
In all optimized Schwarz methods, the matrix $B_i$ cannot be obtained algebraically, as in~\Cref{eq:ai}, due to the boundary conditions. Instead, the local system matrix has to be assembled separately to construct the preconditioner.

\begin{remark} 
	In this paper, we focus on different variants of one-level Schwarz preconditioners. As indicated by the condition number bound in~\Cref{eq:cond:as}, these are not scalable to large numbers of subdomains. When increasing the number of subdomains while keeping
    the computational domain $\Omega$ fixed, $H$ is reduced, and hence, the condition number increases; cf.~\Cref{eq:cond:as}. Even though the convergence of 
    the non-symmetric one-level RAS and ORAS preconditioners cannot be bounded directly via the condition number, they 
    encounter the same scalability issues.
\end{remark}

\section{Coupling deal.II and FROSch}
\label{sec:dealii_and_frosch}
This section discusses the interface between {\dealii} and {\frosch}, which is based on the newly developed {\tpetra}-based {\dealii} 
interface to {\trilinos}.

\subsection{FROSch domain decomposition preconditioner package: Background}

{\frosch} (Fast and Robust Overlapping Schwarz)~\cite{heinlein_frosch_2020,heinlein_parallel_2016-1} is a parallel domain decomposition
preconditioning package, which is part of the {\trilinos} software library~\cite{heroux_overview_2005}. {\frosch} is based on the Schwarz
framework~\cite{toselli_domain_2005} and allows for the construction of Schwarz operators by combining elementary Schwarz operators. A
Schwarz operator is of the form of a preconditioner; for instance, the Schwarz operator corresponding
to the AS preconditioner in~\Cref{eq:as} reads 
\begin{equation} \label{eq:as_op}
	Q_{\rm AS} = M_{\rm AS}^{-1} A.
\end{equation}
Multiple Schwarz operators $Q_1, \ldots, Q_M$ can be combined, for instance, in an additive way,
$$
	Q_{\rm ad} = \sum_{i=1}^{M} Q_i,
$$
or multiplicative way,
$$
	Q_{\rm mult} = I - \prod_{i=1}^M ( I - Q_i ),
$$
where $I$ denotes the identity matrix;
cf.~\cite{toselli_domain_2005}. In this sense, $Q_{\rm AS}$ can also be written as the sum of Schwarz operators $Q_i = P_i A_i^{-1} R_i$
as follows:
$$
	Q_{\rm AS} = \sum_{i=1}^{N} Q_i.
$$
We notice that $M$ is some generic natural number, while $N$ is the number of subdomains as before.

\begin{remark}
	We notice that Schwarz operators are often denoted by the symbol $P_\star$ instead of $Q_\star$; see, for
    instance,~\cite{toselli_domain_2005}. However, we decided on this notation since we already denote the prolongation operators 
    in~\Cref{sec:schwarz} as $P_i$. 
\end{remark}

Another important algorithmic component of {\frosch} are extension-based coarse spaces, which allow for the construction of
preconditioners using only algebraic information or only little geometric information. As discussed in~\Cref{sec:as} and~\Cref{remark:algebraic}, 
the first level of Schwarz preconditioners can be constructed algebraically using only the fully assembled system matrix. For more 
details on the algebraic construction of {\frosch} preconditioners,  see~\cite{heinlein_parallel_2016-1,heinlein_parallel_2016-2,heinlein_fully_2021}.

The main idea of the extension-based coarse spaces is that they do not require an explicit coarse triangulation. 
Still, the coarse spaces 
are constructed from non-overlapping domain
decomposition. In particular, they employ ideas similar to non-overlapping domain decomposition methods, such
as finite element tearing and interconnecting - dual primal (FETI-DP)~\cite{farhat_scalable_2000,farhat_feti-dp_2001} and balancing 
domain decomposition by constraints (BDDC)~\cite{cros_preconditioner_2003,dohrmann_preconditioner_2003}. In particular, the coarse spaces
in {\frosch} are based on generalized Dryja--Smith--Widlund (GDSW)~\cite{dohrmann_family_2008,dohrmann_domain_2008} coarse spaces. These
have been extended to block systems~\cite{heinlein_monolithic_2019}, highly-heterogeneous problems~\cite{heinlein_adaptive_2019},
and composite Discontinuous Galerkin discretizations of multicompartment reaction-diffusion problems~\cite{huynh_gdsw_2024}. Moreover,
reduced-dimension variants~\cite{dohrmann_design_2017,heinlein_improving_2018} and multilevel extensions~\cite{heinlein_multilevel_2023} 
have been developed to reduce the computational cost of solving the coarse problem. Using these techniques, {\frosch} preconditioners have
been
demonstrated to scale up to more than $220,\,000$ cores on the Theta supercomputer; cf.~\cite{heinlein_parallel_2022}. In this
work, we focus more on the one-level preconditioners and the first level of two-level preconditioners.

{\frosch} preconditioners are accessible via the unified solver interface {\stratimikos} of {\trilinos}. They can be constructed from a
fully assembled, parallel distributed sparse matrix in both parallel linear algebra frameworks of Trilinos, {\epetra} or {\tpetra}; 
therefore, the lightweight interface class {\xpetra} is employed. By default, {\frosch} constructs the overlapping domain decomposition
based on the (distributed-memory) parallel distribution of the input matrix using its sparsity pattern; cf.~\Cref{remark:algebraic}. As a 
result, we obtain a one-to-one correspondence of MPI ranks and subdomains. Furthermore, {\frosch} preconditioners accept a parameter list 
as a second input, which can contain additional input parameters. The {\epetra} linear algebra package will soon be deprecated from 
{\trilinos}, so all packages will solely depend on {\tpetra}. Through {\tpetra}, {\frosch} is also able to make use of the 
performance portability model {\kokkos}~\cite{9485033} and the corresponding kernel library 
{\kokkoskernels}~\cite{9485033,rajamanickam_kokkos_2021}, allowing for node-level parallelization on CPUs and GPUs.
The acceleration of {\frosch} preconditioners using GPUs has been demonstrated in~\cite{yamazaki_experimental_2023}.
Through the {\stratimikos} solver interface, {\frosch} preconditioners can easily be combined with iterative Kyrlov solvers from the {\belos} 
package~\cite{bavier_amesos2_2012}. Alternatively, {\frosch} preconditioners and {\belos} solvers can also be called and combined directly 
without using {\stratimikos}. 
For more details on the implementation in {\frosch}, we refer 
to~\cite{heinlein_parallel_2016-1,heinlein_parallel_2016-2,heinlein_monolithic_2019,heinlein_frosch_2020}.

\subsection{\frosch: Our new developments}
  Until now, {\frosch} has provided one- and two-level (R)AS preconditioners, focusing solely on an algebraic approach. 
  Our aim is to extend {\frosch} also to provide an O(R)AS preconditioner. As previously discussed, assembling 
  the local subdomain matrices $B_i$ from~\Cref{eq:oras} requires geometric information.

  Three new classes were added to {\frosch} for the implementation of the O(R)AS preconditioner. The \textit{GeometricOverlappingOperator}
  is derived from the \textit{OverlappingOperator} and acts as the basic building block for either the \textit{GeometricOneLevelPreconditioner} 
  or the \textit{GeometricTwoLevelPreconditioner}. The basic design is similar to the \textit{AlgebraicOverlappingOperator}, but instead of a 
  completely algebraic approach, the \textit{GeometricOverlappingOperator} provides the tools to compute the overlapping domains 
  geometrically, where the overlap is computed based on the dual graph.
  The \textit{GeometricOneLevelPreconditioner} and \textit{GeometricTwoLevelPreconditioner} act similarly to the existing 
  \textit{OneLevelPreconditioner} and \textit{TwoLevelPreconditioner} but are mandatory to expose the tools to compute the overlapping domains to 
  the user. A detailed description of the functions in these three new classes is available as a doxygen 
  documentation\footnote{\url{https://github.com/kinnewig/Trilinos/tree/OptimizedSchwarz}}.

  To construct the ORAS preconditioner, we begin with a parallel distributed triangulation, which is described by lists of vertices, cells, and auxiliary information. Here, each cell is described by $2\cdot dim$ vertices. The auxiliary list has the same length as the cell list and contains additional information, such as material constants. This information must be provided by the finite element program in use. 
  Each rank only stores all vertices, cells, and auxiliary data relevant to itself. Also, the dual graph of how the elements are connected to each other has to be provided from the finite element program. 
  The parallel distributed triangulation is assumed to correspond to a non-overlapping domain decomposition.
  Based on the dual graph, {\frosch} generates the three lists of vertices, cells, and auxiliary information. These
  belong to the overlapping domain decomposition and are returned to the finite element program. Subsequently, 
  local triangulations are created, each corresponding to one subdomain  
  based on the overlapping vertex, cell, and auxiliary data. 
  The local subsystem 
  matrices $B_i$ are assembled with the desired interface conditions on each subdomain. {\frosch} then constructs the O(R)AS 
  preconditioner based on the subsystem matrices $B_i$; cf.~\Cref{eq:oas,eq:oras}.

To summarize, the overlapping domain decomposition is computed automatically based on the user-provided triangulation, 
However, the user must create local triangulations based on the three lists of vertices, cells, and auxiliary information
and, based on that, the local mass matrix and local Neumann matrices with the interface conditions. 
We provide an example of how this can be done in {\dealii}. The preconditioner is automatically computed by {\frosch} based on the local system and interface matrices.
Optionally, the interface conditions could be replaced with different absorbing interface conditions such as PML.

\subsection{\dealii}
  {\dealii} (differential equations analysis library)~\cite{Software:dealii:95,deal2020} is an open-source and widely adopted C++ library specifically designed for  
  solving Partial Differential Equations (PDEs). 
  Its versatility has led to its extensive use in 
  academic research and commercial projects.
  Own examples from ourselves include phase-field fracture propagation in porous media in which up to five PDEs and variational inequalities couple using 
  physics-based discretizations and physics-based solvers~\cite{WheWiLee20}, with the open-source pfm-cracks module~\cite{HeiWi20}, as well as nonstationary fluid-structure interaction optimal control~\cite{WiWo21}, with the open-source optimization toolkit DOpElib~\cite{DOpElib} based on {\dealii}. In fact, repeated computations, such as in optimal control or parameter estimation, are a key reason why the performance of solvers is needed. We also refer to the {\dealii} tutorial page\footnote{\url{https://www.dealii.org/current/doxygen/deal.II/Tutorial.html}} in which numerous further applications are described.
  One of the key attributes of {\dealii} is its high accessibility. This is largely attributed to its comprehensive 
  documentation, active community, and an extensive list of the previously mentioned tutorial steps. These tutorials offer detailed, 
  step-by-step explanations for solving a wide array of equations and problems.  
  The ten main features of {\dealii} are conveniently provided on the main 
  webpage\footnote{\url{https://www.dealii.org/about.html}}.
  
  The construction of one-level O(R)AS preconditioners 
  requires knowledge of which matrix entries belong to the degrees of freedom (dofs) at the interface and, thus, 
  cannot be performed algebraically
  as that information requires knowledge of the grid and how the dofs are distributed on the grid.
  For the two-level O(R)AS preconditioner, we may construct the second level algebraically, as in~\cite{heinlein_parallel_2016-1,heinlein_parallel_2016-2,heinlein_fully_2021}. 
  But since mesh information is readily available in the {\dealii}--{\frosch} interface and required for the construction of the first level, 
  we also provide the index sets of the interface degrees of freedom to {\frosch} as an input. This is because the quality of the algebraic 
  interface reconstruction may depend on the mesh structure and finite elements employed; cf.~\cite{heinlein_fully_2021}. Providing the 
  interface to {\frosch} as input ensures that the exact interface can be employed in the construction of the preconditioners.

\section{Applications}
\label{sec:applications}
In this section, we introduce two challenging applications in order to study the 
performance of our 
algorithms and the newly proposed software implementations.

\subsection{Fluid-Structure Interaction}
Fluid-structure interaction (FSI)~\cite{Ri17_fsi,BaKeTe13,GaRa10} is one of the most challenging PDE systems for the following reasons. It involves the coupling of two physics, namely 
fluid flow and a solid. Large deformations are assumed, which require the modeling of both 
problems in their natural coordinate systems: fluid flow is modeled in Eulerian coordinates, and 
solids are modeled with Lagrangian coordinates. Moreover, the dynamics 
occur on a common interface, which needs to be resolved in the numerical approximation with 
sufficient accuracy 
to prescribe the dynamics to the respective other problem correctly.
In order to couple both systems, different methods exist, one of which is the 
so-called arbitrary Lagrangian-Eulerian (ALE) approach~ 
(\cite{DoFaGi77,HuLiZi81,DoHuePoRo04,FoNo99}) in which fluid flow is re-written 
in general coordinates such that on the interface (and in the solid), Lagrangian
coordinates are utilized while going away; we smoothly transition to Eulerian 
coordinates. 
To this end, we utilize $\hOmega$ as the reference configuration
where the computations occur.
To obtain the physics, the ALE transformation 
$\hcalA$ must be applied to obtain the physical FSI solution in $\Omega$. 
This is also important when the FSI solution is displayed via a visualization program.

In FSI, three solution variables are coupled:
vector-valued velocities $v$,  
scalar-valued pressure $p$, and vector-valued displacements $u$. The resulting discrete 
nonlinear equation systems must first be addressed with some nonlinear solver (often 
Newton's method), then handled internally with a linear solver. The linear systems have a 
saddle-point structure. The solution remains costly despite several solver developments, 
e.g.,~\cite{He04,gee2011truly,CrDeFouQua11,DeparisFortiGrandperrinQuarteroni:2016,
            FortiQuarteroniDeparis:2017,JoLaWi19_fsi,deparis_comparison_2015,
            heinlein_parallel_2016-1,heinlein_parallel_2016,barker_two-level_2010,
            wu_parallel_2011}. 
Specifically, the efficient solution for the three-dimensional case remains a challenge. For this 
reason, FSI has been chosen as one application for our {\dealii}-{\frosch} framework.

We denote by $\Omega\coloneqq\Omega (t) \subset \mathbb{R}^d$, the domain of the FSI 
problem. The domain consists of two time-dependent subdomains $\Omega_f (t)$ and 
$\Omega_s (t)$. The FSI-interface between $\Omega_f (t)$ and $\Omega_s (t)$  is denoted 
by $\Gamma_I (t) = \overline{\partial\Omega_f} (t) \cap \overline{\partial\Omega_s} (t)$.  
The initial (or later reference) domains are denoted by $\hOmega, \hOmega_f$, and 
$\hOmega_s$, respectively, with the interface $\hGamma_i = \overline{\partial\hOmega_f} 
\cap \overline{\partial\hOmega_s}$. Furthermore, we denote the outer boundary by 
$\partial\hOmega = \hGamma = \hGamma_{\text{in}} \cup \hGamma_D \cup \hGamma_{\text{out}}$.
Function values in Eulerian and Lagrangian coordinates are identified by
\begin{equation*}
  u_f(x,t) =: \hu_f(\hx,t) , \quad \text{with } x = \hcalA(\hx,t).
\end{equation*}
The ALE map is defined by
\begin{equation*}
  \hcalA(\hx,t) \coloneqq x = \hx + \hu(x,t)
\end{equation*}
and typically obtained by solving an additional PDE. In this work, this is a nonlinear harmonic model, 
i.e., given $\hat u_s$, at each $t_m$ for $m=1,\ldots,M$, find 
$\hat u_f\colon\bar\hOmega_f\to\mathbb{R}^d$ such that
\begin{equation}
  \label{eq:cauchy_stress_harmonic}
  -\hnabla\cdot(\hat\sigma_{mesh} ) = 0,\quad
  \hat u_f=\hat u_s\text{ on }\hat\Gamma_I,\quad
  \hat u_f=0\text{ on }\partial\hat\Omega_f\setminus\hat\Gamma_I,
\end{equation}
with $\hat\sigma_{mesh} = \hat\alpha_u \hat\nabla \hat u_f $ and $\hat\alpha_u > 0$.
The deformation gradient is given by
\begin{equation*}
  \hF \coloneqq \hnabla \hcalA(\hx,t) = I_d + \hnabla\hu,
\end{equation*}
with $I_d$ being the identity matrix, and its determinant is given by
\begin{equation*}
  \hJ \coloneqq \det(\hF).
\end{equation*}

The following problem statement is a composition of the incompressible Navier-Stokes 
equations, written in ALE coordinates, coupled to geometrically nonlinear elasticity. 
For more details on the derivation, we refer to~\cite{Ri17_fsi}, while our notation 
is based on~\cite{Wi11,Wi13_fsi_with_deal}. For the weak formulation, 
the function space is denoted by $\hX$. Dirichlet boundary conditions 
are built into $\hX$ as usually done, while Neumann conditions 
arise through integration by parts in the weak form. Neumann conditions on the FSI-interface will cancel out as we have the equilibrium of normal stresses. For details, 
we refer to the previously mentioned references. Moreover, the specific 
boundary conditions, such as inflow, outflow, and no-slip conditions, are 
provided in each respective numerical test below. To this end, we have:

\begin{problem}[Variational-monolithic space-time ALE FSI in $\hOmega$]
\label{eq:fsi:ale:harmonic:space_time} 
Let non-homogeneous Dirichlet inflow data, the remaining boundary conditions, and some 
initial data be given. Find a global vector-valued velocity, vector-valued displacements,
and a scalar-valued fluid pressure, i.e., $\hU \coloneqq (\hv,\hu,\hp) \in \hX$ such that:\\
Balance of fluid/solid momentum:
\begin{equation*}
  \begin{aligned}
    \int_I ~ & \Bigg( \left(\hJ\hrho_f  \partial_t \hv
               ,\hpsi^v\right)_{\hOmega_f} 
             + \left(\hrho_f \hJ  (\hF^{-1}(\hv-\hw_A)\cdot\hnabla) \hv,
                     \hpsi^v\right)_{\hOmega_f} 
             + \left(\hJ\hsigma_f\hF^{-T},\hnabla\hpsi^v\right)_{\hOmega_f} \\      
           & - \left\langle \hrho_f \nu_f \hJ(\hF^{-T}\hnabla\hv^T\hn_f)\hF^{-T}, 
                           \hpsi^v \right\rangle_{\hGamma_{\text{out}}} 
             + \left(\hrho_s \partial_t \hv,\hpsi^v\right)_{\hOmega_s}  
             + \left(\hF\hSigma ,\hnabla\hpsi^v\right)_{\hOmega_s} \Bigg)\, dt\\
         + & \left(\hJ\hrho_f(\hv(0) - \hv_{0}), \hpsi^v(0)\right)_{\hOmega_f}
                   + \hrho_s\left(\hv(0) - \hv_{0}, \hpsi^v(0)\right)_{\hOmega_s}
             = 0
  \end{aligned}
\end{equation*}
Mesh/solid 2nd equation:
\begin{equation*}
  \begin{aligned}
    \int_I \Bigg( \left(\hsigma_{\text{mesh}} ,\hnabla\hpsi^u\right)_{\hOmega_f}
      + 
      \hrho_s \left(\partial_t\hu-\hv,\hpsi^u\right)_{\hOmega_s} \Bigg) \, dt
      + \hrho_s \left(\hu(0) - \hu_{0}, \hpsi^u(0)\right)_{{\hOmega_s}} 
      = 0
  \end{aligned}
\end{equation*}
Fluid mass conservation:
\begin{equation}
  \begin{aligned}
      \int_I \Bigg( \left(\hnabla\cdot(\hJ\hF^{-1}\hv),\hpsi_f^p\right)_{\hOmega_f} \Bigg)\, dt
       = 0
  \end{aligned}
\end{equation}
for all $\hPsi = (\hpsi^v,\hpsi^u,\hpsi^p) \in \hX^0$. 
The fluid constitutive stress tensor is given by 
\begin{equation*}
  \hsigma_f \coloneqq -\hp_f I_d +\hrho_f\nu_f(\hnabla \hv_f \hF^{-1}
  + \hF^{-T}\hnabla \hv_f^T),
\end{equation*}
modeling a Newtonian fluid.
The viscosity and the density of the fluid are denoted by $\nu_f$ and $\hrho_f$, respectively. 
The solid constitutive stress tensors are given by the second Piola-Kirchhoff stress tensor
\begin{equation}
  \hSigma = \hSigma_s(\hu_s) = 2\mu \hE + \lambda tr(\hE) I_d, 
\end{equation}
where the Green-Lagrange strain tensor is given by 
\begin{equation*}
  \hE = \frac{1}{2} (\hF^T \hF - I_d).
\end{equation*}
Here, $\mu$ and $\lambda$ are the Lam\'e coefficients for the solid. The solid density is
denoted by $\hrho_s$. The resulting solid material model is of St.~Venant Kirchhoff type.
\end{problem}
The Rothe method (horizontal methods of lines) is employed 
for the discretization: first, discretizing in time, then 
in space. For temporal discretization, 
finite difference schemes are used, namely the (implicit) 
backward Euler scheme (BE), which is strongly A-stable but only first-order 
and dissipative. A second choice is the Crank-Nicolson scheme, which is A-stable, 
second-order, less dissipative, and well-suited for energy conservation problems. 
The time-discretized equations are the starting point for the 
Galerkin discretization in space. To this end, we construct finite dimensional
subspaces $\hat X_h^0 \subset \hat X^0$ to find an approximate solution to the 
continuous problem. In the context of monolithic ALE formulations, the 
computations are done on the reference configuration $\hOmega$. We use two- or 
three-dimensional shape-regular meshes. A mesh consists of quadrilateral or 
hexahedral elements $\hat K$. They form a non-overlapping cover of the 
computation domain $\hat\Omega\subset\mathbb{R}^d$, $d=2,3$. The corresponding 
mesh is given by $\mathcal{\hat T}_h = \{ \hat K \}$. The discretization parameter 
in the reference configuration is denoted by $\hat h$ and is a cell-wise constant
that is given by the diameter $\hat h_{\hat K}$ of the cell $\hat K$. The 
resulting nonlinear discrete systems are solved with a line-search-assisted 
Newton's method; see~\cite{Wi11,Wi13_fsi_with_deal}. 
The resulting linearized equation systems are solved using our new
{\dealii}-{\frosch} 
software interface. Lastly, we mention that, in~\cite{Wi13_fsi_with_deal}, our {\dealii}
FSI code was published and later moved to Github\footnote{\url{https://github.com/tommeswick/fsi}}, 
serving as a starting point for our 2D FSI implementation. The 3D FSI code can be 
implemented accordingly, with its prior results shown in~\cite{Wi11}[Section 4.4].

\subsection{Time-Harmonic Maxwell's Equations} 
\label{sec:maxwell}
Maxwell's equations~\cite{Bk:Fey:63, Bk:Monk:2003} are utilized in many modern research 
fields, from magnetic induction tomography (MIT) in healthcare~\cite{Art:Zolgharni:MIT:2009}, 
geo-electromagnetic modeling in geophysics~\cite{Art:Grayver:Geoelectromagnetic:2015}, 
quantum computing~\cite{Art:Kues:23} and optics/photonics~\cite{Art:Melchert:Soliton:2023}.

\begin{problem}[Time-Harmonic Maxwell's Equations]
    Let us consider $\Omega \subset \mathbb{R}^3$, a bounded Lipschitz domain with a sufficiently smooth 
    boundary $\Gamma = \Gamma^{\text{inc}} \cup \Gamma^{\infty}$. On $\Gamma^{\infty}$, an absorbing 
    boundary condition is prescribed.
    One possible choice is the Robin boundary condition
    \begin{equation}
      \label{eq:robin-boundary}
      \mu^{-1}\left(\vec{n} \times \left(\operatorname{curl}\left(\vec{u}\right)\right)\right) 
      - \mathsf{i} \sqrt{\varepsilon} \omega \left(\vec{n} \times \left(\vec{u} \times \vec{n}\right)\right)
      = 0,
    \end{equation}
    which corresponds to the first-order approximation of the Sommerfeld radiation conditions.
    On $\Gamma^{\text{inc}}$, a boundary condition for a given 
    incident electric field is prescribed.. The goal is to find the electric field 
    ${\vec u} \in \mathbf{H}_\text{curl}(\Omega)$ such that, for all 
    $\vec \varphi \in \mathbf{H}_\text{curl}(\Omega)$, the following equation holds:
    \begin{equation}
      \label{eq:THM_weak}
      \begin{split}
        \int_\Omega \left( \mu^{-1} \operatorname{curl} \left( {\vec u} \right) \cdot \operatorname{curl} \left( \vec \varphi \right)
        - \varepsilon \omega^2 {\vec u} \cdot\vec \varphi \right) ~ \mathsf{d} x
        + \mathsf{i} \sqrt{\varepsilon} \omega \int_{\Gamma} 
        (\vec{n}\times (\vec{u} \times \vec{n}))\cdot (\vec{n}\times (\vec{\varphi}\times\vec{n} )) ~ \mathsf{d} s, \\
        = \int_{\Gamma^\text{inc}} 
        (\vec{n}\times (\vec{u}^\text{inc} \times \vec{n})) \cdot (\vec{n}\times (\vec{\varphi}\times \vec{n} )) ~ \mathsf{d} s. 
      \end{split}
    \end{equation}
    Here, ${\vec u}^\text{inc}$ with $ \vec{n}\times {\vec u}^\text{inc}\in L^2(\Gamma^{\text{inc}},
    \mathbb{C}^d)$ is a given incident electric field, $\mu \in \mathbb{R}^+$ is the relative magnetic 
    permeability, $\varepsilon \in \mathbb{R}^+$ is the relative 
    permittivity, $\omega = \frac{2 \pi}{\lambda}$ is the wavenumber, and $\lambda \in \mathbb{R}^{+}$ 
    is the wavelength. Finally, $\mathsf{i}$ is the imaginary number.
    \end{problem}
\begin{remark}
    System~\Cref{eq:THM_weak} is called time-harmonic because the time dependence can be expressed by $e^{\mathsf{i} \omega \tau}$, where $\tau > 0$ represents the time. We avoid the need for a time-stepping method, where possibly a high number of
    time steps are required until the harmonic state is reached. We refer the reader to~\cite{Bk:Monk:2003} for a detailed derivation of the time-harmonic Maxwell's equations.
\end{remark}
Solving the indefinite time-harmonic Maxwell's equations is 
challenging, as the linear system resulting from a Finite 
Element (FE) discretization is ill-posed
as the corresponding system matrix has positive and negative eigenvalues.
One encounters similar difficulties in the numerical solution of the ill-posed Maxwell's equations, which also occur in the Helmholtz equation~\cite{Art:Gander:Helmholtz:2012}.
Therefore, specialized techniques have to be employed. 
DDMs yield good performance compared with other techniques, especially when 
optimized interface conditions are used; cf.~\cite{Art:DoGaGe:09}. When dealing with Maxwell's equations 
in the context of FE, it is necessary to use N\'ed\'elec elements to preserve the curl-conforming behavior of the 
underlying function spaces. As N\'ed\'elec elements are oriented, one has to take care of the orientations 
of each edge and face when creating the (overlapping) subdomains to match the 
original orientation. For a description of the orientation problem in the context of N\'ed\'elec elements, 
we refer the reader to~\cite{Art:Kinnewig:23}.
All these technical details must be considered while constructing an O(R)AS preconditioner for Maxwell's 
equations, making this a challenging application.

\section{Numerical studies} 
\label{sec:results}
In this section, we first demonstrate the performance and ease of use of the {\dealii}-{\frosch}
interface by solving different FSI problems with the use of RAS. Secondly, we demonstrate 
the O(R)AS extension of the {\dealii}-{\frosch} interface on the example of time-harmonic Maxwell's 
problems. To solve the linear problems on each subdomain, the direct solver \texttt{MUMPS}~\cite{Software:MUMPS}
was used. All results presented in this section were computed on a server with a 2 $\times$ AMD EPYC 7H12 64-Core processor with 1024\,GB of RAM.

\subsection{Fluid-Structure Interaction} \label{sec:fsi_results}
First, we consider a series of FSI problems using the two-level RAS preconditioner 
provided by {\frosch} in {\dealii}. 
Due to the equal-order discretization with $Q_1\times Q_1\times Q_1$ elements, we can employ a monolithic two-level Schwarz preconditioner similar to~\cite{heinlein_frosch_2022}. 
However, we use a two-level RAS preconditioner with a simple second level in which the coarse basis functions have constant nodal values within each partition of a unique partitioning of the finite element nodes; cf.~\cite{nicolaides_deflation_1987}. This coarse space is algebraic, and it can be directly constructed using, for instance, the $\hat P_i$ prolongation operators of the RAS preconditioner~\Cref{eq:pou1}, or the corresponding binary scaling matrices $D_i$; cf.~\Cref{eq:pou2}. Even though this coarse space might not be optimal, we
chose this two-level preconditioner as it still performed better than a one-level preconditioner, 
especially for a higher number of subdomains. 
For a brief comparison between the one-level and the two-level method, see the solver analysis in~\Cref{subsec:fsi-1}.

\subsubsection{Example 1: 2D FSI-1 benchmark}
\label{subsec:fsi-1}
This benchmark was proposed in~\cite{HrTu06b}, and results from various groups 
are shown in~\cite{BuSc10}. The characteristic feature of this benchmark is 
a steady-state solution, which we will approximate via a backward Euler time-stepping scheme.

\begin{figure}
  \begin{center}
    \includegraphics[width=0.95\textwidth]{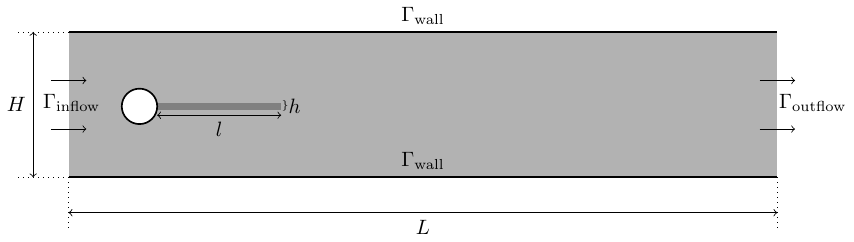}
  \end{center}
  \caption{
    \label{fig:fsi-domain}
    Flow around an obstacle to which an elastic beam is attached. 
    In the two-dimensional case, the domain is a rectangle, and the obstacle is modeled as a 
    cycle with the center $C(0.20,0.20)$,
    and the lower right corner of the elastic beam is point $A(0.60,0.19)$. 
    In the three-dimensional case, the domain is a rectangular cuboid with depth $D$. The 
    obstacle is modeled as a cylinder with the center $C$, and the front lower right corner of 
    the elastic beam is point $A$.
  }
\end{figure}

\paragraph{Configuration, boundary conditions, quantities of interest}
The computational domain is defined by a rectangle with dimensions $L=2.50\;\unit{m}$ and 
$H=0.41\;\unit{m}$; see~\Cref{fig:fsi-domain}. Within this domain, there is an obstacle that is 
slightly offset from the center, represented by the circle with the center 
$C=(0.20\;\unit{m},0.20\;\unit{m})$ and the radius $r=0.05\;\unit{m}$. On the circle, an elastic 
beam of length $l=0.35\;\unit{m}$ and height $h=0.02\;\unit{m}$ is attached, and the lower 
right end of that beam is located A = $(0.60\;\unit{m}, 0.19\;\unit{m})$.
As material parameters, we consider a fluid density of $\rho_f = 1.00 \cdot 10^3\;\unit{\sfrac{kg}{m^3}}$ 
the viscosity $\nu_f=1.00 \cdot 10^{-3}\;\unit{\sfrac{m^2}{s}}$. We consider a density 
of $\rho_s = 1.00 \cdot 10^3 \;\unit{\sfrac{kg}{m^3}}$ for the solid and a Lam\'e coefficient 
of $\mu = 0.50 \cdot 10^6\;\unit{\sfrac{kg}{m~s^2}}$.
For the boundary conditions on $\Gamma_\text{inflow}$, we consider 
$v_\text{inflow} = 0.20\;\unit{\sfrac{m}{s}}$ to be the inflow velocity.
On the outflow boundary, we have the do-nothing condition~\cite{HeRaTu96}, namely, normal stresses equal 
to zero, which results in the correction term 
$\langle \hrho_f \nu_f \hJ(\hF^{-T}\hnabla\hv^T\hn_f)\hF^{-T}, 
                         \hpsi^v \rangle_{\hGamma_{\text{out}}}$
due to the symmetric stress tensor. Lastly, we prescribe 
no-slip conditions (homogeneous Dirichlet) on the top and bottom walls 
of the channel, i.e., $v = 0$.
We use $\Delta t = 1.00\;\unit{s}$ for the time step size.

\begin{figure}
  \begin{center}
    \includegraphics[width=0.95\textwidth]{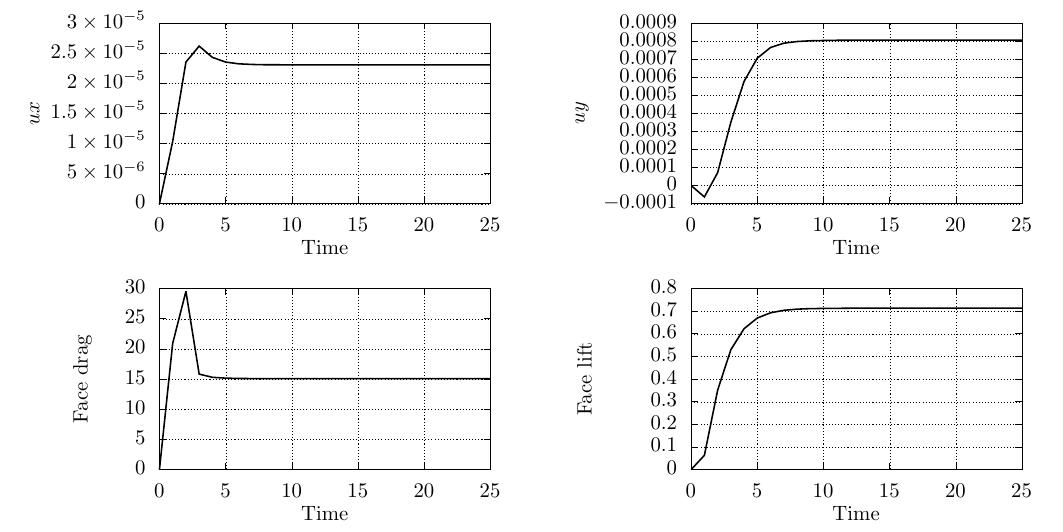}
  \end{center}
  \caption{
    \label{fig:FSI-1-results}
    FSI-1 results.
    \textbf{Top left:} 
      $x$-displacement in $A(0.6,0.19)$.
    \textbf{Top right:} 
      $y$-displacement in $A(0.6,0.19)$.
    \textbf{Bottom left:}
      Face drag measured around the cylinder and the elastic beam.
    \textbf{Bottom right:}
      Face lift measured around the cylinder and the elastic beam.
    }
\end{figure}

\paragraph{Physics results}
To verify our results, we consider four quantities of interest~\cite{HrTu06b}, namely, the position $A:=A(0.60,0.19)$; see~\Cref{fig:fsi-domain}. This is the lower right edge of the elastic beam, where we 
plot the $x$ and $y$ components individually. Moreover, we consider the values for the Face lift and the Face drag. 
These are evaluated around the cylinder and the elastic beam. For the mathematical formulae, we refer 
to~\cite{HrTu06b,Wi13_fsi_with_deal}. The results are visualized in~\Cref{fig:FSI-1-results}.

\paragraph{Solver analysis}
To evaluate the performance of the two-level RAS preconditioner, we examine the number of GMRES iterations necessary to 
solve the problem when we modify the number of subdomains and keep the number of dofs per subdomain constant. 
In two dimensions, this is achieved by a single global refinement and a fourfold increase in the number of 
subdomains. 
As we use a partition of unity with one function per subdomain for the coarse space, the coarse space size is given by 
$(2 \cdot dim + 1)(ranks)$ because one constant value is assumed for each variable on any subdomain. 

It is worth noting that the grid is first partitioned by \texttt{p4est} and then refined.
Secondly, we consider the number of GMRES iterations dependent on the overlap between the 
subdomains. The results are shown in~\Cref{tab:iteration-count-fsi-1}. 

For the two-level method and for an overlap of $\sfrac{\delta}{H} = 20\%$ and 128 subdomains 
on average 
468 GMRES iterations were required to solve the linear system, 
where a maximum of 589 GMRES iterations was encountered, and it took around $47.50\;\text{minutes}$. 
In comparison, the one-level method took an average of 
555
iterations for the same case, 
where a maximum of 856 GMRES iterations was encountered, and it took around $56.78\;\text{minutes}$. 

\begin{figure}
  \centering
  \includegraphics[width=0.95\textwidth]{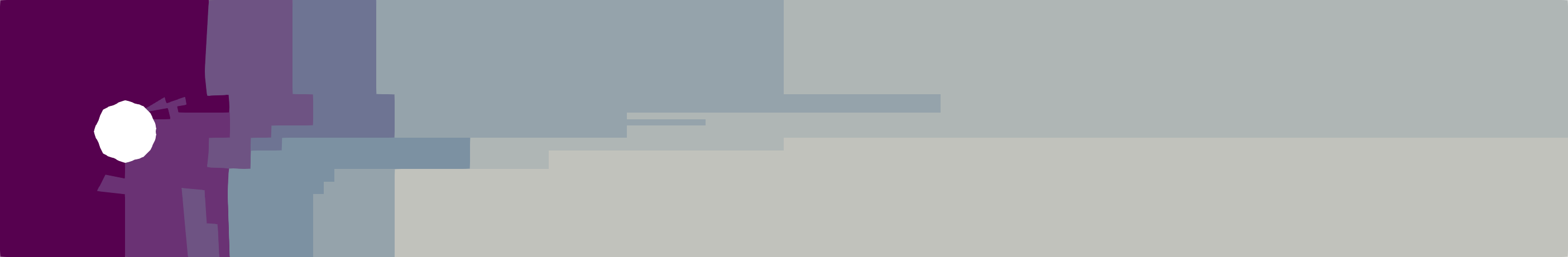}
  \caption{
    \label{fig:domain-decomposition-fsi}
    The domain decomposition of the FSI-1 benchmark, as obtained by applying \texttt{p4est}. 
    There are non-connected subdomains, as well as very long subdomains, which is less than optimal
    for a domain decomposition.
  }
\end{figure}

Ideally, with a two-level RAS preconditioner and a large number of subdomains, the number of GMRES iterations 
would remain constant. However, given our second level’s simplicity, we observe an increase in the number of 
GMRES iterations, roughly proportional or slightly more, each time we quadruple the number of subdomains.
This increase can be attributed to two factors. First, the default grid partitioning tool from {\dealii} is 
\texttt{p4est}\cite{Software:p4est:11}, 
which was not designed with the application to domain decomposition 
in mind, leading to less-than-optimal partitioning
(see~\Cref{fig:domain-decomposition-fsi}).

Second, the large time step size seems to contribute 
significantly to the high GMRES iteration counts; see also~\Cref{tab:timestepsize} for iteration counts for 
varying time step sizes as well as the results in~\Cref{subsec:fsi-2}.
When the overlap is doubled, the number of GMRES iterations should be cut in half for large enough numbers of
subdomains. Even though we observe a strong correlation between the overlap and the number of GMRES iterations, 
we do not observe this behavior. We suspect this might be related to the mesh partitioning obtained from~\texttt{p4est}.

\begin{table}
  \begin{center}
    \begin{tabular}{rrllllll} \toprule
                 & \multicolumn{3}{c}{Overlap}                            \\
      ranks & dofs    & $\sfrac{\delta}{H}=5\,\%$ & $t^{\text{walltime}}_{5\%}$ & $\sfrac{\delta}{H}=10\,\%$ & $t^{\text{walltime}}_{10\%}$ & $\sfrac{\delta}{H}=20\,\%$ & $t^{\text{walltime}}_{20\%}$\\ \midrule
      2          & $  4\,070$  & $19.94~(25) $   & $1.0\;\unit{min}$ & $13.81  ~(16) $ & $1.0\;\unit{min}$ & $10.03  ~(12) $  & $ 1.1\;\unit{min}$\\
      8          & $ 15\,620$  & $91.15~(107)$   & $1.6\;\unit{min}$ & $59.27  ~(74) $ & $1.7\;\unit{min}$ & $39.01  ~(50) $  & $ 2.4\;\unit{min}$\\
      32         & $ 61\,160$  & $331.15~(549)$  & $3.9\;\unit{min}$ & $171.77 ~(202)$ & $3.7\;\unit{min}$ & $102.70 ~(118)$  & $ 5.4\;\unit{min}$\\
      128        & $242\,000$  & $>1000        $ & $-$               & $>1000        $ & $-$               & $468.30 ~(589)$  & $47.5\;\unit{min}$\\ \bottomrule
    \end{tabular}
  \end{center}
  \caption{
    \label{tab:iteration-count-fsi-1}
    Overview of the wall times and the average number of GMRES iteration steps required to solve the first 25 time steps
    of the FSI-1 problem in dependence on the number of used ranks, which corresponds to the number 
    of subdomains and the overlap between the subdomains. The number in the brackets shows the highest 
    number of GMRES iterations encountered.
    We consider the average of GMRES iterations here as a Newton method, and a time-step method is 
    performed to solve the FSI problem.    
  }
\end{table}

\begin{table}
  \centering
  \begin{tabular}{rlllll} \toprule
    time step size   & 1.0000 & 0.5000 & 0.2500 & 0.1250 & 0.0625 \\ \midrule
    GMRES iterations & 346.63 & 255.16 & 221.62 & 199.60 & 178.99 \\ \bottomrule
  \end{tabular}
  \caption{
    \label{tab:timestepsize}
    The number of GMRES iterations in dependency of the time step size, on the example of the FSI-1 benchmark with $61\,160$ dofs, 
    $\sfrac{\delta}{H}=5\%$ overlap and 128 subdomains.
  }
\end{table}

\subsubsection{Example 2: 2D FSI-2 benchmark}
This benchmark was also proposed in~\cite{HrTu06b}, but yields periodic solutions and the von Karman vortex street.
Numerous results for comparison were published in~\cite{HrTu06b} and~\cite{Wi11,Wi13_fsi_with_deal}.

\paragraph{Configuration, boundary conditions, quantities of interest}
We consider the same domain as for the FSI-1 example in~\Cref{subsec:fsi-1}; 
see~\Cref{fig:fsi-domain}. 
As material parameters, we consider the same parameters for the fluid as in the FSI-1 example. 
For the solid, we consider a density  of $\rho_s = 1.00 \cdot 10^4\;\unit{\sfrac{kg}{m^3}}$ 
and a Lam\'e coefficient of $\mu = 0.50 \cdot 10^6\;\unit{\sfrac{kg}{m~s^2}}$.
For the boundary conditions on $\Gamma_\text{inflow}$, we consider 
$v_\text{inflow} = 1.00\;\unit{\sfrac{m}{s}}$ to be the inflow velocity.
The other boundary conditions are the same as above.
For the time step size we use $\Delta t = 1.00\cdot10^{-2}\;\unit{s}$.

\begin{figure}
  \centering
  \includegraphics[width=0.95\textwidth]{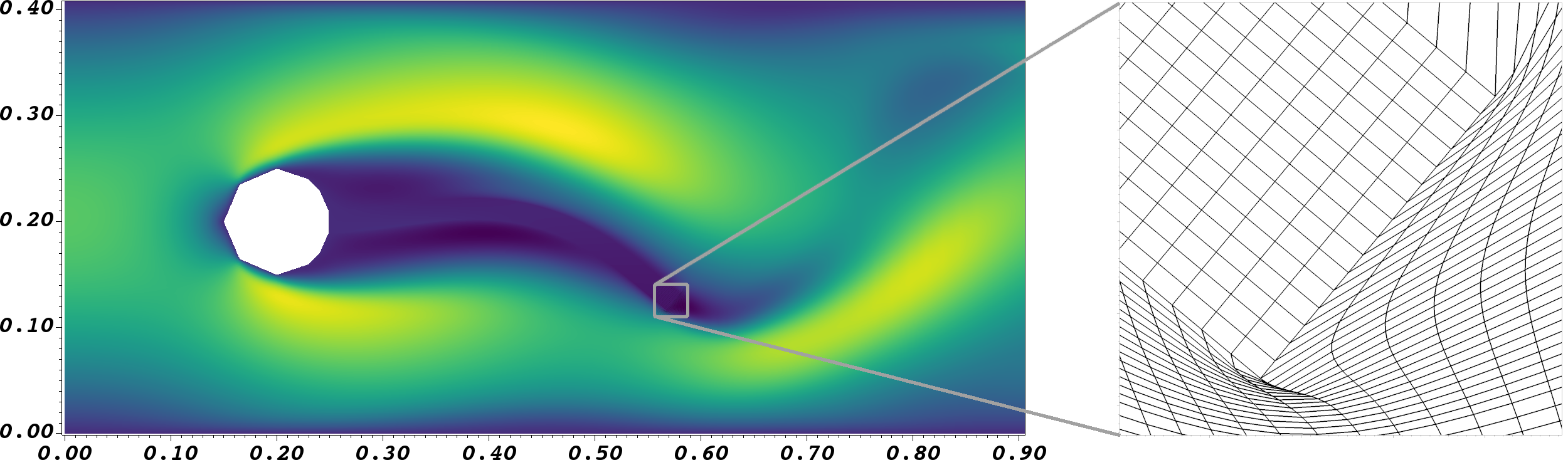}
  \caption{
    \label{fig:FSI-2-xvelocity}
      \textbf{Left:} Plot of the resulting velocity into $x$-direction from the flow around the cylinder at time step $t=10\;\unit{s}$.
      \textbf{Right:} The mesh deformation at the tip of the elastic beam at time step $t=10\;\unit{s}.$ 
  }
\end{figure}

\begin{figure}
  \begin{center}
    \includegraphics[width=0.95\textwidth]{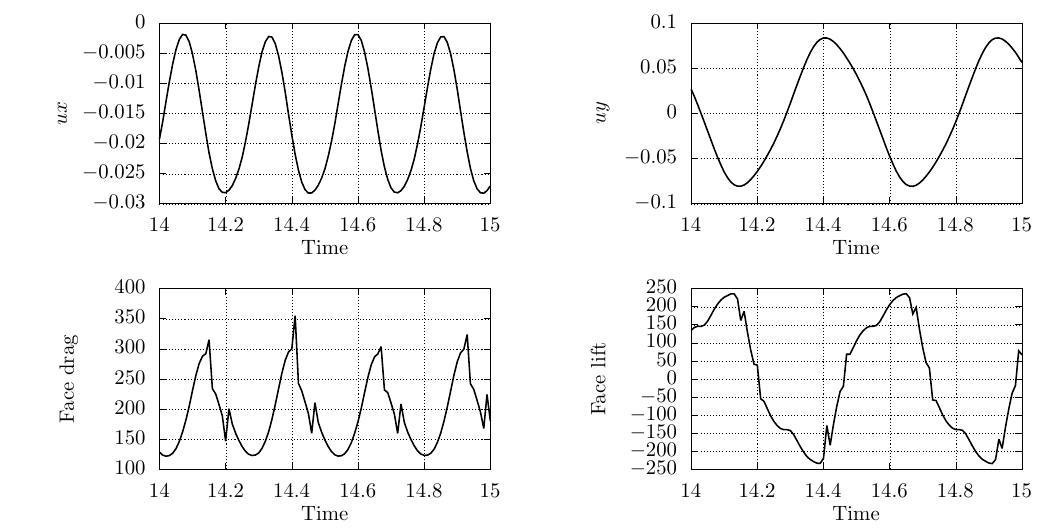}
    \includegraphics[width=0.95\textwidth]{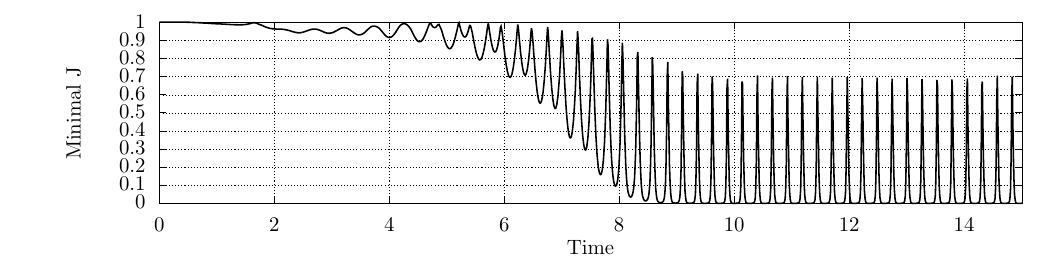}
  \end{center}
  \caption{
    \label{fig:FSI-2-results}
    FSI-2 results. 
    \textbf{Top left:} 
      $x$-displacement in $A(0.6,0.19)$.
    \textbf{Top right:} 
      $y$-displacement in $A(0.6,0.19)$.
    \textbf{Mid left:}
      Face drag.
    \textbf{Mid right:}
      Face lift.
    \textbf{Bottom:}
      Time evolution of the minimal $\hJ$ (determinant of the deformation gradient $\hF$) over all $1500$ time steps.
  }
\end{figure}

\paragraph{Physics results}
In this example, a periodic state is obtained~\cite{HrTu06b, Wi13_fsi_with_deal}. The FSI-2 benchmark is the most 
difficult of all three proposed benchmarks as large deformations of the elastic beam occur, requiring a robust 
mesh motion model to avoid mesh distortion~\cite{Wi11}; a graphical illustration of the highest deflection 
is given in~\Cref{fig:FSI-2-xvelocity}. To represent this periodic behavior, we again focus on the point 
$A(0.6,0.19)$, that is, the lower left point of the beam.  
Our analysis focuses on the periodic state, and all findings are compared to~\cite{HrTu06b,Wi13_fsi_with_deal,Wi11}. Therefore, we only consider the time $t>14\;\unit{s}$ 
as proposed in~\cite{HrTu06b}. In~\Cref{fig:FSI-2-results}, the periodic motion in terms of the $y$ displacements can be observed. 
Our findings for drag and lift slightly differ from the literature values because we use equal-order 
finite elements and a simple pressure stabilization, which enters into the stress evaluation and thus 
slightly alters drag and lift values. Overall, the behavior of all four quantities of interest is in an 
acceptable range and confirms the correctness of the implemented FSI problem.
Finally, the minimal $\hJ$ is plotted as well, and the nominal minimal value is 
$8\cdot 10^{-4}$,
which shows that the chosen mesh motion model still works but is close to its limit.

\paragraph{Solver analysis}
As above, we study the correlation between the number of subdomains and GMRES iterations and the
correlation between overlap and GMRES iterations. The results are shown in~\Cref{tab:iteration-count-fsi-2}.
Here, we consider the average values over the first $200$ time steps. 
Here, the number of GMRES iterations increases less than for the FSI-1 example. 
Going from $32$ to $128$ subdomains, the number of GMRES iterations doubles, even though the number 
of subdomains quadruples. 
This is because of the smaller time step size. As shown in~\Cref{tab:timestepsize}, there is a strong correlation between the time step size and the number of GMRES iterations.

\begin{table}
  \begin{center}
    \begin{tabular}{rrllllll} \toprule
                 & \multicolumn{3}{c}{Overlap}                             \\
      ranks & dofs    & $\sfrac{\delta}{H}=5\,\%$ & $t^{\text{walltime}}_{5\%}$ & $\sfrac{\delta}{H}=10\,\%$ & $t^{\text{walltime}}_{10\%}$ & $\sfrac{\delta}{H}=20\,\%$ & $t^{\text{walltime}}_{20\%}$\\ \midrule
      2          & $  4\,070$  & $19.99  ~(22) $ & $6.5\;\unit{min}$  & $13.9 ~(15) $ & $6.60\;\unit{min}$  & $11.0  ~(12) $ & $7.1\;\unit{min}$  \\
      8          & $ 15\,620$  & $68.32  ~(82) $ & $10.0\;\unit{min}$ & $47.6 ~(54) $ & $10.77\;\unit{min}$ & $33.7  ~(39) $ & $15.3\;\unit{min}$ \\
      32         & $ 61\,160$  & $129.98 ~(172)$ & $15.0\;\unit{min}$ & $88.1 ~(112)$ & $17.57\;\unit{min}$ & $63.1  ~(75) $ & $27.5\;\unit{min}$ \\
      128        & $242\,000$  & $240.42 ~(355)$ & $77.2\;\unit{min}$ & $163.4~(211)$ & $91.48\;\unit{min}$ & $117.1 ~(155)$ & $128.7\;\unit{min}$ \\ \bottomrule
    \end{tabular}
  \end{center}
  \caption{
    \label{tab:iteration-count-fsi-2}
      Overview of the wall times and the average number of GMRES iteration steps required to solve the first 200 time steps
      of the FSI-2 problem in dependence on the number of used ranks and the overlap between the subdomains. 
      The number in the brackets shows the highest number of GMRES iterations encountered.
  }
\end{table}

The computation of 1500 time steps with $242\,000$ dofs took roughly $36.39\;\text{hours}$ with $20\%$ overlap.
The average number of Newton steps per time step is 
$5$.
The minimal encountered number is $1$, and the highest encountered number 
is $16$ (where the mesh is close to degeneration with $\hJ$ close to zero). As a result, on average, 
the solution and assembly of each linear system took around $16.67\;\text{seconds}$.

\subsubsection{Example 3: 2D elastic lid-driven cavity}
\label{subsec:fsi-2}
In this example, we consider a test case in which the solid covers a larger part of the overall domain~\cite{Du06}. 
We adapt the setting to have the solid cover a larger part of the domain.

\begin{figure}[H]
  \begin{center}
    \includegraphics[width=0.65\textwidth]{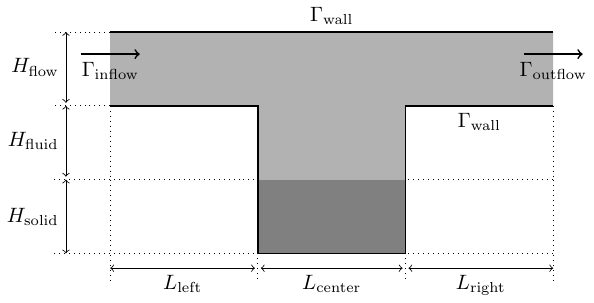}
  \end{center}
  \caption{
    \label{fig:lit-driven-cavity}
      Configuration of the lid-driven cavity with the elastic bottom (dark gray) in the cavity.
   }
\end{figure}

\paragraph{Configuration, boundary conditions, quantities of interest}
For the computational domain, see~\Cref{fig:lit-driven-cavity}. Here, we consider the length
$L_\text{left} = L_\text{center} = L_\text{right} = 0.50\;\unit{m}$ and the height 
$H_\text{solid} = H_\text{fluid} = H_\text{flow} = 0.25\;\unit{m}$.
As a material parameter, we consider the fluid to have the same parameters as in FSI-1; see~\Cref{subsec:fsi-1}. 
We consider a density structure of $\rho_s = 1.00\;\unit{\sfrac{kg}{m^3}}$ for the solid. 
The Lam\'e coefficient is $\mu = 1.00 \cdot 10^2\;\unit{\sfrac{kg}{m~s^2}}$.
For the boundary conditions on $\Gamma_\text{inflow}$, we consider 
$v_\text{inflow} = 2.00\;\unit{\sfrac{m}{s}}$ to be the inflow velocity.
The other boundary conditions are the same as in the previous examples:
do-nothing outflow and no slip on the walls.
We use $\Delta t = 1.00\;\unit{s}$ for the time step size.

\begin{figure}[H]
  \begin{center}
    \includegraphics[width=0.95\textwidth]{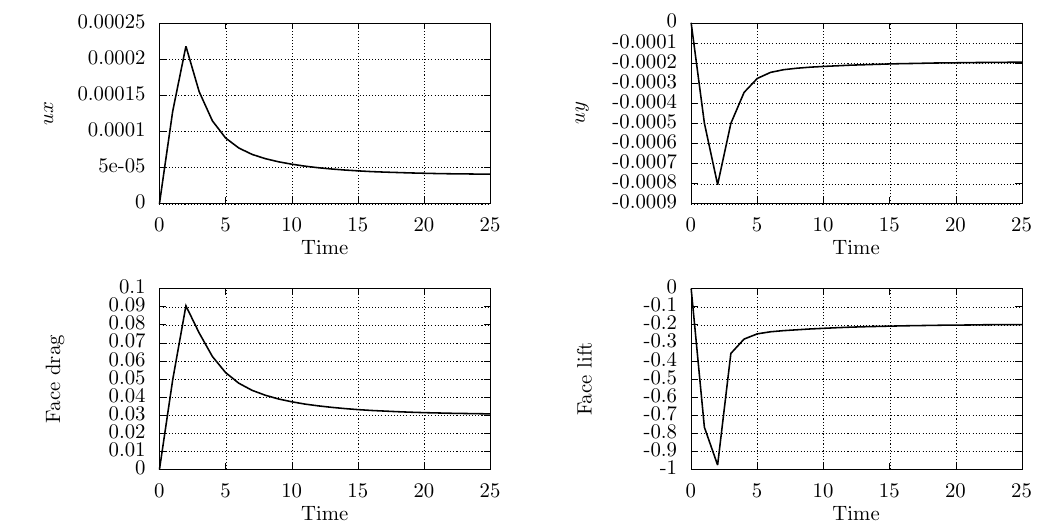}
  \end{center}
  \caption{
    \label{fig:FSI-cavity-results}
    Results from the lid-driven cavity with elastic bottom. 
    \textbf{Top left:} 
      $x$-displacement in $A = (0.75, 0.25)$.
    \textbf{Top right:} 
      $y$-displacement in $A = (0.75, 0.25)$.
    \textbf{Bottom left:}
      Face drag (normal stresses in $x$-direction on the interface).
    \textbf{Bottom right:}
      Face lift (normal stresses in $y$-direction on the interface).
  }
\end{figure}

\paragraph{Physics results}
In this case, a vortex forms in the cavity, deforming the elastic solid at the bottom of the solid.
To verify our results, we consider the point $A = (0.75\;\unit{m}, 0.25\;\unit{m})$, which sits 
in the center of the interface between the solid and the liquid; see~\Cref{fig:lit-driven-cavity}.
In~\Cref{fig:FSI-cavity-results}, the $x$ and $y$ components of the point $A$ are plotted. 
Moreover, the normal stresses (named Face drag and Face lift) are shown.
The results correspond to what we expected, as one can observe that the point $A$ is shifted down 
by the vortex forming in the liquid above the solid.

\paragraph{Solver analysis}
  In this example, a larger part of the domain is covered by the solid, which is beneficial for the RAS preconditioner.
  As we can see, when we compare the results in~\Cref{tab:iteration-count-fsi-cavity} to the 
  above examples, the number of GMRES iterations grows less in dependency on the 
  number of subdomains than it does in the FSI-1 example, even though a time step size of $\Delta t = 1.00\;\unit{s}$ 
  was used here as well. 

\begin{table}
  \begin{center}
    \begin{tabular}{rrllllll} \toprule
                 & \multicolumn{3}{c}{Overlap}                             \\
      ranks & dofs    & $\sfrac{\delta}{H}=5\,\%$ & $t^{\text{walltime}}_{5\%}$ & $\sfrac{\delta}{H}=10\,\%$ & $t^{\text{walltime}}_{10\%}$ & $\sfrac{\delta}{H}=20\,\%$ & $t^{\text{walltime}}_{20\%}$\\ \midrule
      2          & $  6\,885$  & $ 31.17  ~(39)$ & $28.13\;\unit{min}$ & $ 17.31  ~(19)$ & $28.10\;\unit{min}$ & $ 15.03  ~(16)$ & $28.45\;\unit{min}$ \\
      8          & $ 26\,565$  & $103.71 ~(115)$ & $27.38\;\unit{min}$ & $ 42.74  ~(50)$ & $27.13\;\unit{min}$ & $ 30.74  ~(35)$ & $27.78\;\unit{min}$ \\
      32         & $104\,325$  & $223.13 ~(244)$ & $27.48\;\unit{min}$ & $101.56 ~(110)$ & $26.73\;\unit{min}$ & $ 71.18  ~(76)$ & $28.10\;\unit{min}$ \\
      128        & $413\,445$  & $472.00 ~(555)$ & $52.82\;\unit{min}$ & $252.76 ~(287)$ & $42.75\;\unit{min}$ & $170.85 ~(189)$ & $48.33\;\unit{min}$ \\ \bottomrule
    \end{tabular}
  \end{center}
  \caption{
    \label{tab:iteration-count-fsi-cavity}
    Overview of the wall times and the average number of GMRES iteration steps required to solve the first 50 time steps
    of the lid-driven-cavity problem in dependence on the number of used ranks and the overlap between the subdomains. 
    The number in the brackets shows the highest number of GMRES iterations encountered.
  }
\end{table}

  In the computation with $20\%$ overlap, the average number of Newton steps per time point $t_n$ is 
  $2$. The minimal encountered number 
  of Newton iterations is $1$, and the highest encountered number is $4$. As a result, on average, the solution and assembly of each linear 
  system took around $33.72\;\text{seconds}$.

\subsubsection{Example 4. 3D-FSI: 3D elastic bar behind a square cross-section}
\begin{figure}
  \centering
  \includegraphics[width=0.95\textwidth]{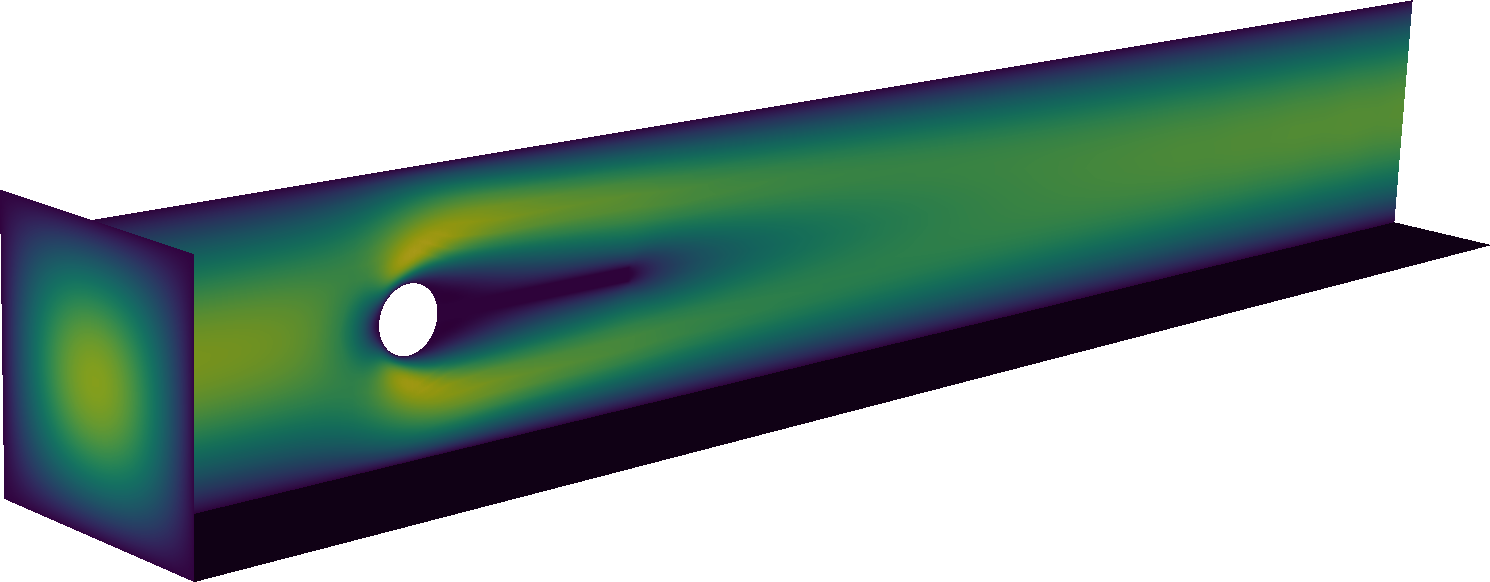}
  \caption{
    \label{fig:FSI-3d-xvelocity}
      3D-FSI example. The plot of the $x$-velocity at timestep $t=25\;\unit{s}$ displaying the final steady-state solution.
  }
\end{figure}

This 3D-FSI example is inspired by~\cite{Wi11}, which is an extension of the FSI-1 benchmark into three spatial dimensions. Specifically, a steady-state solution is obtained; see~\Cref{fig:FSI-3d-xvelocity}. 
A cross-section of the computational domain is displayed 
in~\Cref{fig:fsi-domain}.

\paragraph{Configuration, boundary conditions, quantities of interest}
For the dimensions of the cuboid, the height is $H=0.41\;\unit{m}$, the length $L=2.80\;\unit{m}$, 
and the depth is $D=0.41\;\unit{m}$. The center of the cylinder is $C=(0.50\;\unit{m}, 0.20\;\unit{m}, 0.00\;\unit{m})$
with the radius $r=0.05\;\unit{m}$.
The front, lower right edge of the elastic beam is located at $A=(0.90\;\unit{m}, 0.19\;\unit{m}, 0.10\;\unit{m})$.
The other end is attached to the cylinder. The elastic beam has a length of $l=0.45\;\unit{m}$, a height of 
$h=0.02\;\unit{m}$, and a depth of $d=0.20\;\unit{m}$.
For the material parameters and boundary conditions, we assume the same as in the previous example 
FSI-1 in~\Cref{subsec:fsi-1}.
We use $\Delta t = 1.00\;\unit{s}$ for the time step size.

\begin{figure}
  \begin{center}
    \includegraphics[width=0.95\textwidth]{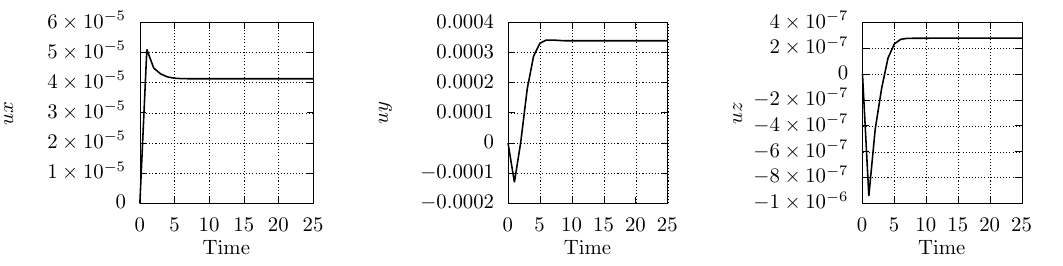}\\
    \includegraphics[width=0.95\textwidth]{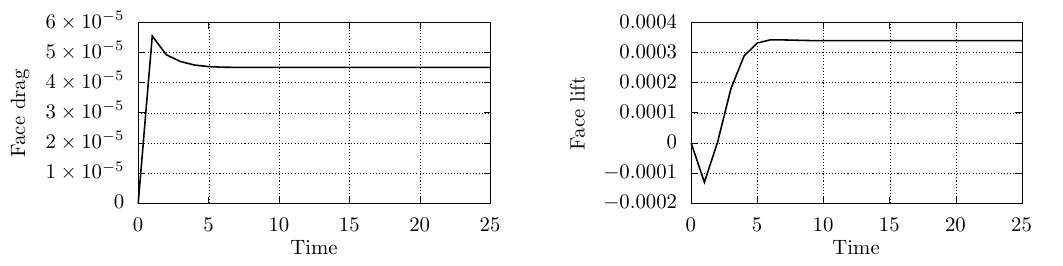}
  \end{center}
  \caption
  {
    \label{fig:FSI-3d-results}
    Results from the 3D-FSI problem. 
    \textbf{Top left:} 
      $x$-displacement of the point A of the elastic beam, i.e., the front lower right, a static state is assumed after some time.
    \textbf{Top center:} 
      $y$-displacement of the point A of the elastic beam.
    \textbf{Top right:} 
      $z$-displacement of the point A of the elastic beam.
    \textbf{Bottom left:}
      Global Face drag plotted against the time.
    \textbf{Bottom right:}
      Global Face lift plotted against the time.
  }
\end{figure}

\paragraph{Physics results}
After some initialization time, the system assumes a steady state, as in the FSI-1 example. For the evaluation, we consider the $x$, $y$, and $z$ components of point A, as well as
the Face drag and the Face lift in~\Cref{fig:FSI-3d-results}.
In the time-dependent position of point $A$, the Face drag, and the Face lift, we can observe that the system converges against a steady state, as expected. The velocity in the 
$x$-direction inside the cuboid is shown in~\Cref{fig:FSI-3d-xvelocity}.

\paragraph{Solver analysis}
We study the correlation between the number of subdomains and the number of GMRES iterations for the three-dimensional case, 
similar to the previous examples. To keep the number of dofs per subdomain constant, we globally refine once and multiply 
the number of subdomains by a factor of eight. 
As before, we also study the correlation between the number of GMRES 
iterations and the size of the overlap between subdomains. If we increase the number of subdomains by a factor of eight, the 
number of GMRES iterations grows by a factor of four or less. The number of GMRES iterations grows less in dependence on 
the number of subdomains compared to the example in~\Cref{subsec:fsi-1}, even though the main difference between 
these two examples is the spatial dimension. 
The results are shown in~\Cref{tab:iteration-count-fsi-3d}.

For $20\%$ overlap, the average number of Newton steps per time point $t_n$ is $6$ iterations. The minimal encountered number of Newton 
iterations is $6$, and the highest encountered number is $9$ iterations. This implies that, on average, the solution and assembly of each linear system 
took around $4.41\;\text{minutes}$.

\begin{table}
  \begin{center}
  \begin{tabular}{rrllll} \toprule
               & \multicolumn{2}{c}{Overlap}                             \\
    ranks & dofs & $\sfrac{\delta}{H}=10\,\%$ & $t^{\text{walltime}}_{10\%}$ & $\sfrac{\delta}{H}=20\,\%$ & $t^{\text{walltime}}_{20\%}$\\ \midrule
    2     & $    19\,278$ & $23.83  ~(26) $ & $ 17.4\;\unit{min}$ & $16.21  ~(20) $ & $ 17.4\;\unit{min}$ \\
    16    & $   135\,660$ & $95.13  ~(104)$ & $ 59.5\;\unit{min}$ & $74.31  ~(83) $ & $103.8\;\unit{min}$ \\
    128   & $1\,014\,552$ & $266.22 ~(291)$ & $424.7\;\unit{min}$ & $190.09 ~(209)$ & $687.7\;\unit{min}$ \\ \bottomrule
  \end{tabular}
  \end{center}
  \caption{
    \label{tab:iteration-count-fsi-3d}
    Overview of the wall times and the average number of GMRES iteration steps required to solve the first 25 time steps of the 3D-FSI example
    in dependence on the number of used ranks and the overlap between the subdomains. 
    The number in the brackets shows the highest number of GMRES iterations encountered.
    As we observed poor convergence for $\sfrac{\delta}{H}=5\,\%$, we skip this case here.
  }
\end{table}

\subsection{Time-Harmonic Maxwell's Equations}
We consider Maxwell's equations as a second class of examples to demonstrate the performance of the one-level ORAS preconditioner.
As discussed in \cite{Art:DoGaGe:09}, classical (R)AS preconditioners do not yield optimal performance in these kinds of problems. 
As the boundary condition for the local subdomain matrices $B_i$, we chose an easy-to-implement Robin-type interface condition 
shown in~\Cref{eq:robin-boundary} and the Robin parameter to be $\alpha=1.0$.
To ensure the correctness of the results computed by our Maxwell solver, we conducted a comparative analysis with a reference
Maxwell solver computing exactly the same problem statement. As reference solver, we use the Maxwell solver described in~\cite{Art:Kinnewig:22:DD26}, which employs MUMPS; therefore, 
we could compute reference solutions for the 2D and 3D simple waveguide. In all cases, it holds for both numerical solutions 
$\| u_{\text{ref}} - u_{\text{ORAS}}\|_{L_{\infty}} < 10^{-12}$, where
$u_{\text{ref}}$ is the numerical solution computed by the reference Maxwell solver, and $u_\text{ORAS}$ is the
numerical solution of the Maxwell solver developed in this paper. It can be inferred that 
the physics are correct employing our new solver.

As for the third example, we consider a real-world application closely related to the problem discussed in~\cite{Art:Kinnewig:22:DD26}.
In all experiments, we consider a wavelength of $632.8\;\unit{nm}$, corresponding to the light
emitted by a Helium-Neon laser. As a material with a higher refractive index, we are using quartz (SiO$_2$),
which has, for the considered wavelength, a refractive index of 
$n_{\text{SiO}_2@632.8\;\unit{nm}}=1.4570$ ($\mu_{\text{SiO}} = 1.0000$ and $\varepsilon_{\text{SiO}_2} = n_{\text{SiO}_2}^2$).
As a cladding material, we assume air with the refractive index $n_\text{air}=1.0000$, 
($\mu_{\text{air}} = 1.0000$, $\varepsilon_{\text{air}} = 1.0000$).

In this subsection, we consider a one-level preconditioner; therefore, we expect the number of GMRES iterations to increase with the 
number of subdomains. Even though~\Cref{eq:cond:as} cannot be applied directly to the case of using GMRES for Maxwell's equations, 
the condition number bound for the Laplace would indicate a condition number of $\kappa(M^{-1}_{AS}) \leq C (1 + \sfrac{1}{\beta H^2})$ for
$\sfrac{\delta}{H}\coloneqq \beta$ to be constant. 
Then, if we quadruple the number of subdomains in a two-dimensional setting with a structured domain decomposition, $H$ is cut in half. 
Therefore, the condition number is multiplied by a factor of four, i.e., the number of GMRES iterations will increase as well.

Finally, we mention that we consider N\'ed\'elec elements of the lowest order, and our findings are compared to~\cite{Art:Kinnewig:22:DD26}.
Here, we partition the grid using \texttt{METIS}~\cite{karypis_fast_1998}. 

\subsubsection{Example 1: 2D Simple Waveguide}
\label{subsec:2d-waveguide}
As a first example, we consider a simple model of a two-dimensional waveguide.

\begin{figure}
  \begin{minipage}{0.45\textwidth}
    \begin{center}
      \includegraphics[width=0.80\textwidth]{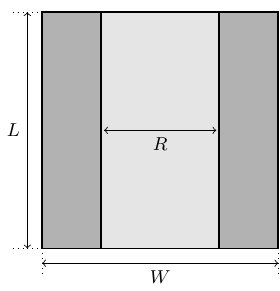}
    \end{center}
  \end{minipage}
  \hfill
  \begin{minipage}{0.45\textwidth}
    \begin{center}
      \includegraphics[width=0.95\textwidth]{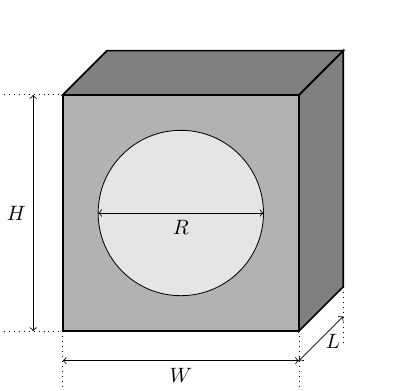}
    \end{center}
  \end{minipage}
  \caption{
    \label{fig:waveguides-geometrie}
    \textbf{Left:} Model of a simplified two-dimensional waveguide, the light gray area 
    indicates the core, i.e., the domain with a higher refractive index, and the dark gray 
    area indicates the cladding, i.e., the domain with a lower refractive index.
    \textbf{Right:} Model of a simplified glass fiber, i.e., a simple waveguide. Like the 
    two-dimensional waveguide, the light gray domain indicates the core and the dark gray area 
    indicates the cladding.
  }
\end{figure}

\paragraph{Configuration, boundary conditions, quantities of interest}
The waveguide is modeled by the domain $\Omega = (0, W) \times (0, L) = (0, 2) \times (0, 6)\; \mu m$, 
which includes a second rectangular subdomain of a material with a higher refractive index
with width $R=0.8\;\unit{\mu m}$ at its center. 
The geometry of the model is shown in~\Cref{fig:waveguides-geometrie} (left). Moreover, 
$\Gamma_{\text{inc}} =  (0, 2) \times \{ 0 \} \; \mu m$ is the boundary with the incident boundary 
condition. All other boundaries, represented by $\Gamma_\infty$, are characterized by absorbing boundary
conditions, for which we choose homogeneous Robin boundary conditions.
The incident electric field is modeled by 
$u_{\text{inc}} = \operatorname{exp}\left(\sfrac{-100}{\mu m^2} (x^2)\right) \vec{e}_x$, 
where $\vec{e}_x$ is the unit vector in the $x$-direction.

\begin{figure}
  \begin{minipage}{0.45\textwidth}
    \begin{center}
      \includegraphics[width=0.95\textwidth]{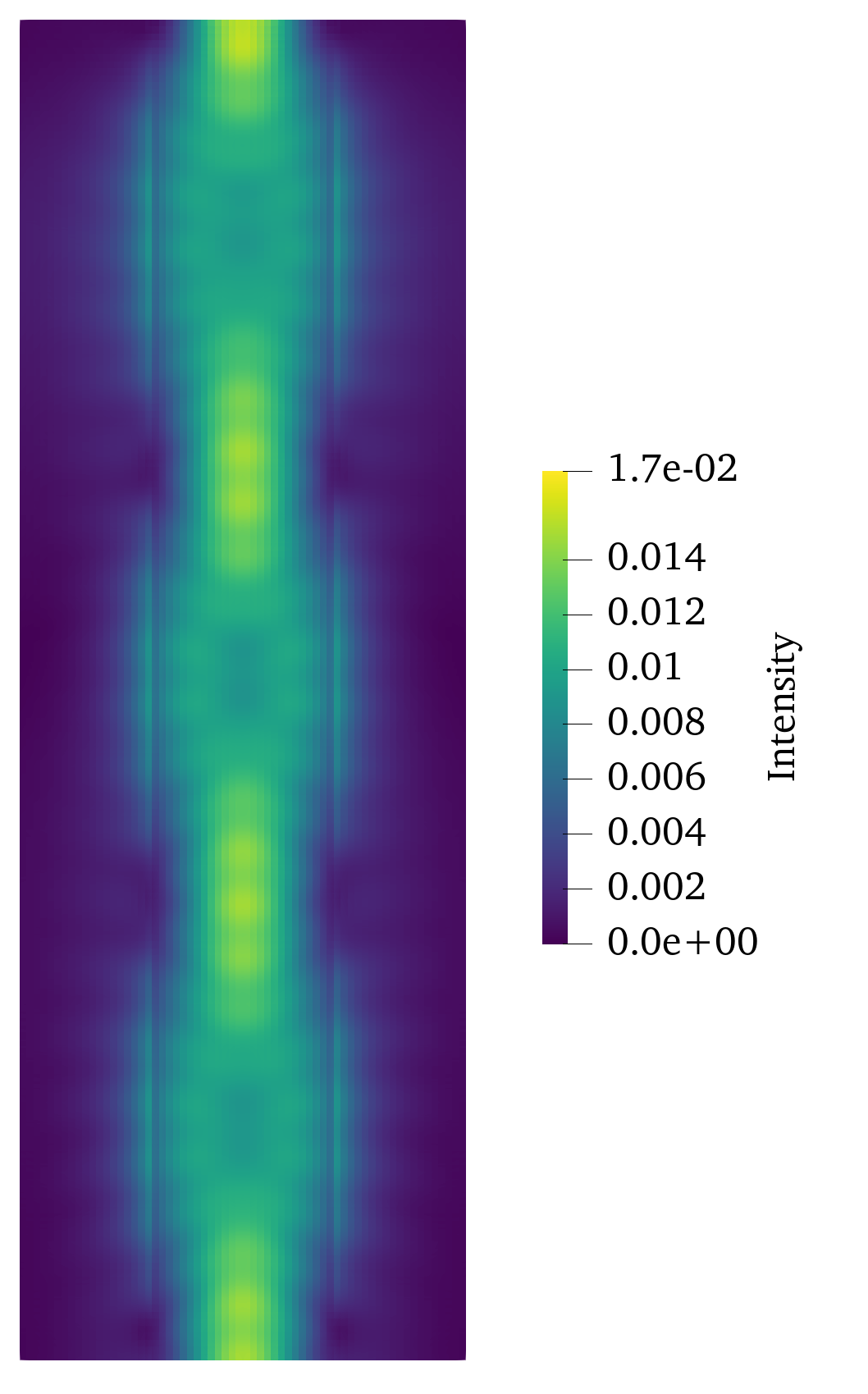}
    \end{center}
  \end{minipage}
  \hfill
  \begin{minipage}{0.45\textwidth}
    \begin{center}
      \includegraphics[width=0.95\textwidth]{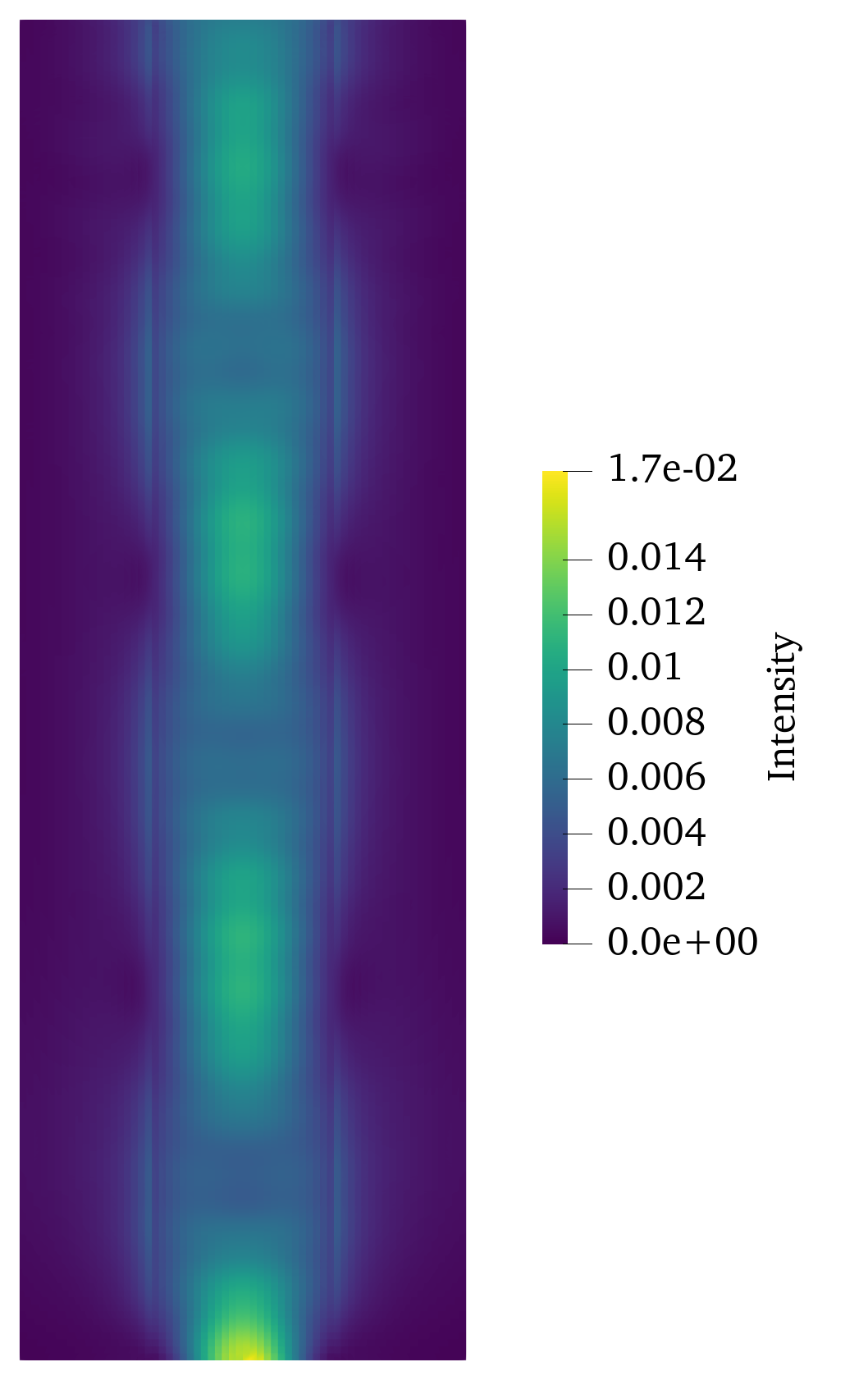}
    \end{center}
  \end{minipage}
  \caption{
    \label{fig:waveguide-solution}
    \textbf{Left:} Intensity plot of the electric field in the two-dimensional waveguide.
    \textbf{Right:} Cross-section of the three-dimensional waveguide at the plane 
    $(0, 2) \times \{1\} \times (0, 6)\;\unit{\mu m}$.
    Intensity plot of the electric field inside of the three-dimensional waveguide.
  }
\end{figure}

\paragraph{Physics results}
  \Cref{fig:waveguide-solution} (left) shows the intensity distribution of the electric field forming inside the two-dimensional waveguide. 
  As we only consider the intensity of the electric field here, we can see the formation of waves inside the waveguide, as one would expect.

\paragraph{Solver analysis}
  Similar to the FSI problems, we are interested in the correlation of GMRES iterations with
  the overlap and the number of subdomains; the corresponding results are shown in~\Cref{tab:itertation-count-maxwell-simplewaveugide-2d}. 
  To keep the number of dofs per 
  subdomain constant, we used one global refinement and increased the number of subdomains by 
  a factor of four. The number of GMRES iterations grows at a 
  similar rate as the number of subdomains, which is the expected behavior for a one-level method.

\begin{table}
  \begin{center}
    \begin{tabular}{rrllllll} \toprule
                 & \multicolumn{6}{c}{Overlap}                            \\
      ranks & dofs    & $\sfrac{\delta}{H}=5\,\%$ & $t^{\text{walltime}}_{5\%}$ & $\sfrac{\delta}{H}=10\,\%$ & $t^{\text{walltime}}_{10\%}$ & $\sfrac{\delta}{H}=20\,\%$ & $t^{\text{walltime}}_{20\%}$\\ \midrule
      2   & $    832$ & $12 $ & $0.08\;\unit{s}$ & $12 $ & $0.09\;\unit{s}$ & $8  $ & $0.09\;\unit{s}$ \\
      8   & $ 3\,200$ & $37 $ & $0.14\;\unit{s}$ & $33 $ & $0.16\;\unit{s}$ & $23 $ & $0.17\;\unit{s}$ \\
      32  & $12\,544$ & $100$ & $0.37\;\unit{s}$ & $86 $ & $0.40\;\unit{s}$ & $70 $ & $0.44\;\unit{s}$ \\
      128 & $49\,664$ & $405$ & $7.73\;\unit{s}$ & $223$ & $4.94\;\unit{s}$ & $176$ & $9.16\;\unit{s}$ \\\bottomrule
    \end{tabular}
  \end{center}
  \caption{
    \label{tab:itertation-count-maxwell-simplewaveugide-2d}
    Overview of the wall times and the number of GMRES iterations required to solve the two-dimensional simple waveguide 
    using the ORAS preconditioner in dependence on the number of used ranks, which corresponds to the 
    number of subdomains and the overlap between the subdomains.
  }
\end{table}

\subsubsection{Example 2: 3D Simple Waveguide}
  This example is the three-dimensional variant of the example discussed above. We consider a simplified 
  glass fiber model consisting of one material with a higher refractive index in the center
  surrounded by a material with a lower refractive index.

\paragraph{Configuration, boundary conditions, quantities of interest}
The waveguide is modeled inside the domain 
$\Omega = (0, W) \times (0, H) \times (0, L) = (0, 2) \times (0, 2) \times (0, 6)\;\unit{\mu m}$, 
which includes a cylinder
of radius $R=0.4\;\unit{\mu m}$ structure at its core. The geometry 
of the model is shown in~\Cref{fig:waveguides-geometrie} (right). Moreover, 
$\Gamma_{\text{inc}} =  (0, 2) \times (0, 2) \times \{ 0 \} \;\unit{\mu m}$ is the boundary 
with the incident boundary condition. All other boundaries, represented by $\Gamma_\infty$, 
are characterized by absorbing conditions. As mentioned above, we employ homogeneous Robin conditions as 
absorbing boundary conditions.
The incident electric field is modeled by 
$u_{\text{inc}} = \operatorname{exp}\left(\sfrac{-100}{\mu m^2} (x^2+y^2)\right) \vec{e}_x$, 
where $\vec{e}_x$ is the unit vector in the $x$-direction.

\paragraph{Physics results}
  The intensity distribution of the electric field forming inside the three-dimensional waveguide is shown in~\Cref{fig:waveguide-solution} (right).

\paragraph{Solver analysis}
As in the previous example, we study the correlation between the number of subdomains and the number of 
GMRES iterations. To keep the number of elements per subdomain constant in the three-dimensional setting, the number 
of subdomains has to be multiplied by a factor of eight for every global refinement.
The results are shown in~\Cref{tab:itertation-count-maxwell-simplewaveugide-3d}.
\begin{table}
  \begin{center}
    \begin{tabular}{rrllllll} \toprule
            & \multicolumn{6}{c}{Overlap}                            \\
      ranks & dofs    & $\sfrac{\delta}{H}=5\,\%$ & $t^{\text{walltime}}_{5\%}$ & $\sfrac{\delta}{H}=10\,\%$ & $t^{\text{walltime}}_{10\%}$ & $\sfrac{\delta}{H}=20\,\%$ & $t^{\text{walltime}}_{20\%}$\\ \midrule
      $2  $ & $81\,056  $   & $12 $ & $ 42.6\;\unit{s}$ & $10 $ &  $43.8\;\unit{s}$ & $9  $ & $ 46.8\;\unit{s}$ \\
      $16 $ & $618\,816 $   & $94 $ & $ 61.6\;\unit{s}$ & $75 $ & $ 66.6\;\unit{s}$ & $48 $ & $ 85.6\;\unit{s}$ \\
      $128$ & $4\,833\,920$ & $397$ & $356.0\;\unit{s}$ & $155$ & $246.0\;\unit{s}$ & $129$ & $374.0\;\unit{s}$ \\\bottomrule
    \end{tabular}
  \end{center}
  \caption{
    \label{tab:itertation-count-maxwell-simplewaveugide-3d}
    Overview of the wall times and the number of GMRES iterations required to solve the three-dimensional simple waveguide 
    using the ORAS preconditioner dependent on the number of used ranks and the overlap between the subdomains.
  }
\end{table}

\subsubsection{Example 3: Y-Beamsplitter}
As a last example, we consider a real-world problem, that is, a so-called y-beam splitter,
where we have one waveguide as input and two waveguides as output. This example is closely 
related to the problem discussed in \cite{Art:Kinnewig:22:DD26}.

\begin{figure}
  \begin{minipage}{0.45\textwidth}
    \begin{center}
      \includegraphics[width=0.95\textwidth]{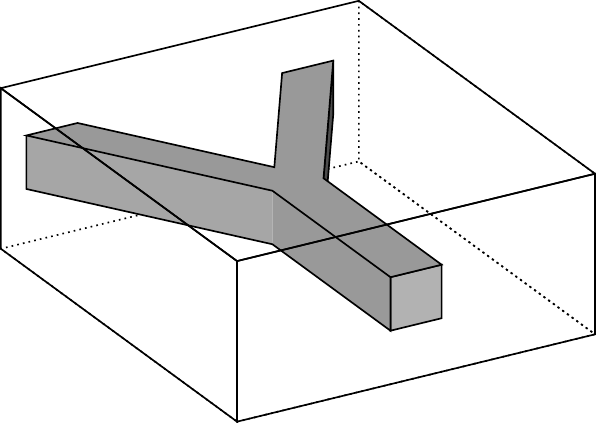}
    \end{center}
  \end{minipage}
  \hfill
  \begin{minipage}{0.45\textwidth}
    \begin{center}
      \includegraphics[width=0.95\textwidth]{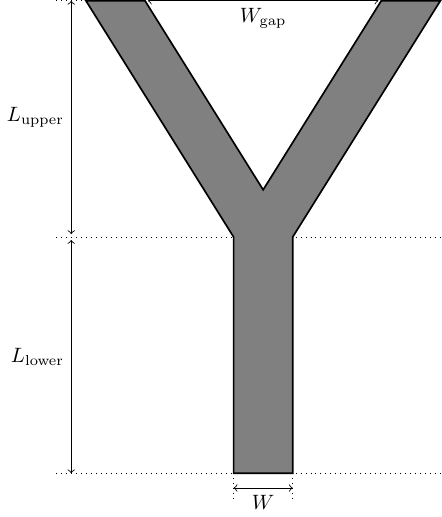}
    \end{center}
  \end{minipage}
  \caption{
    \label{fig:y-beamsplitter-geometry}
      \textbf{Left:} The orientation of the y-beam splitter inside of the cladding.
      \textbf{Right:} A top-down view of the y-beam splitter, with labels for important geometric characteristics.
  }
\end{figure}

\paragraph{Configuration, boundary conditions, quantities of interest}
The y-beam splitter is modeled inside the domain 
$\Omega = (0.0, 3.2) \times (0.0, 1.2) \times (0.0, 6.4)\;\unit{\mu m}$, 
where the geometry of the y-beam splitter itself is shown in~\Cref{fig:y-beamsplitter-geometry} (left).
The width and the height of the waveguides the y-beam splitter consists of are
$W=0.8\;\unit{\mu m}$,
$L_\text{lower} = 0.8\;\unit{\mu m}$, $L_\text{upper} = 5.6;\unit{\mu m}$, 
and the gap between the ends of the outgoing waveguides is $W_\text{gap} = 0.4\;\unit{\mu m}$; see~\Cref{fig:y-beamsplitter-geometry} (right).
The incident electric field is modeled by 
$u_\text{inc} = \operatorname{exp}\left(\sfrac{-100}{\mu m^2} (x^2 + y^2) \right) \vec{e}_x$,
where $\vec{e}_x$ is the unit vector in the $x$-direction.

\begin{figure}
  \begin{center}
    \includegraphics[width=0.95\textwidth]{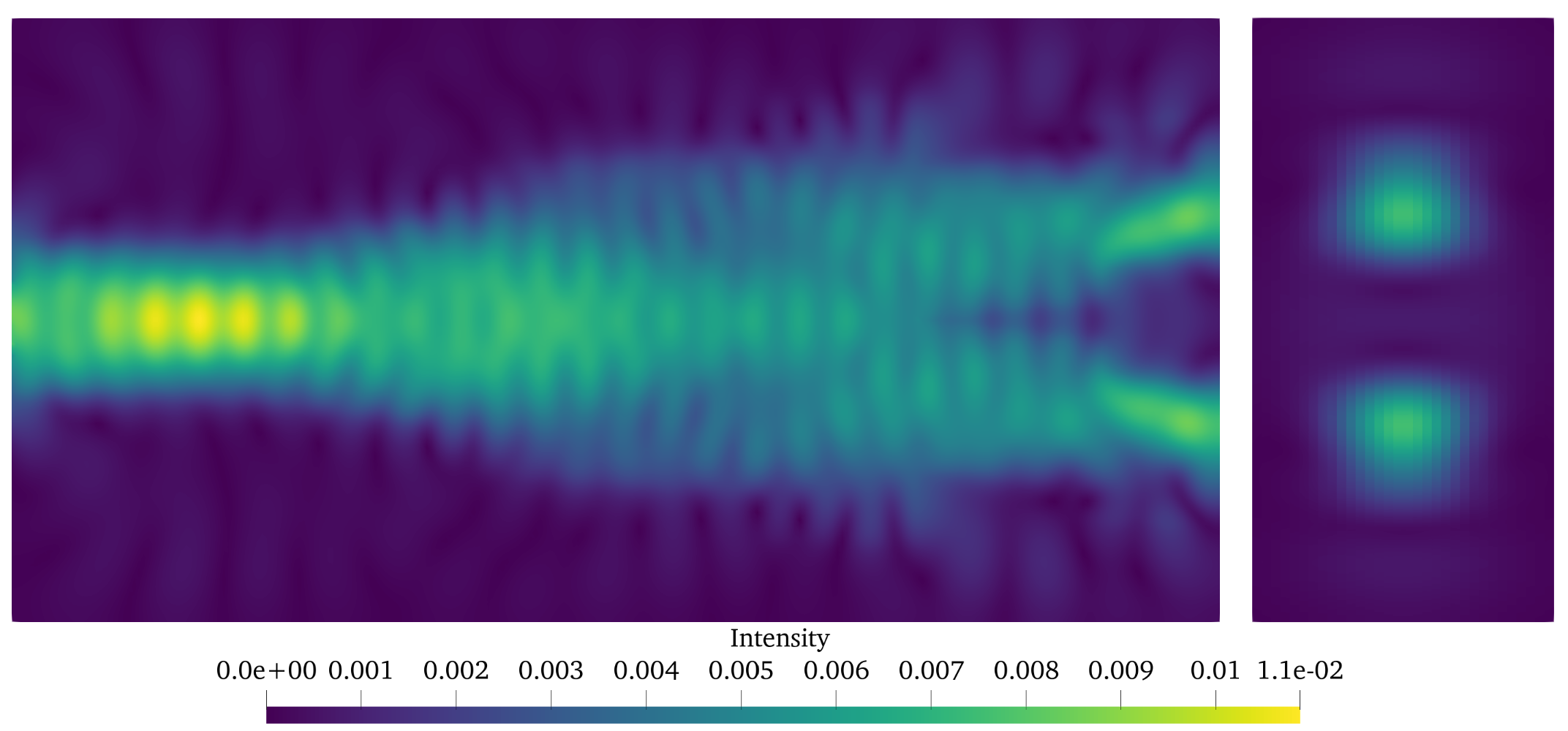}
  \end{center}
  \caption{
    \label{fig:y-beamsplitter-result}
    \textbf{Left:} Cross-section of the y-beam splitter at the plane $(0.0, 3.2) \times \{0.6\} \times (0.0, 6.4)\;\unit{\mu m}$ 
    showing the Intensity distribution of the electric field inside of the y-beam splitter. 
    \textbf{Right:} The intensity distribution of the output after splitting one input into two outputs.
  }
\end{figure}

\paragraph{Physics results}
  The solution of the y-beam splitter is challenging for two reasons. 
  First, a sufficiently small mesh element size must be chosen based on the y-beam splitter size and the electric field's wave-like nature. 
  At least four degrees of freedom per wave are required;
  otherwise, the wave-like nature cannot be resolved.
  Second, the splitting point introduces a lot of 
  reflection, making the problem harder to solve. We can also observe both of these problems here, as discussed
  in the solver analysis below. But with a sufficiently large overlap, i.e., $\delta\geq10\,\%$, the y-beam splitter
  can be efficiently solved with the ORAS preconditioner. 
  Our results are displayed graphically in~\Cref{fig:y-beamsplitter-result}.

\paragraph{Solver analysis}
  In this example, a convergence analysis is more cumbersome. As discussed above, a sufficiently small step width
  is required to capture the wave-like nature of the electric field. The two smaller examples, with $112\,964$,
  and $856\,168$ dofs, respectively, are insufficient to resolve the single waves. 
  Therefore, it would be reasonable to start with an even higher initial refinement, requiring additional computational 
  resources, which will not provide more insight into the performance of the {\dealii} - {\frosch} interface and these are 
  therefore left for future HPC simulation studies of practical applications within the excellence cluster PhoenixD.
  The relatively low initial refinement explains why we observe
  a large jump in the iteration count in the case of an overlap of $5\,\%$ between 16 and 128 subdomains. For $5\,\%$ 
  overlap, we only observe poor convergence. Our results are shown in~\Cref{tab:itertation-count-maxwell-ybeamsplitter}.
  This highlights that a sufficient overlap must be chosen for more complex problems. 
  Alternatively, more advanced interface conditions can be employed, e.g., high-order approximations as 
  proposed in \cite{Art:Bouajaji:15} or PML-based interface conditions \cite{Art:Dolean:TransmissionMaxwell:24}.

\begin{table}
  \begin{center}
    \begin{tabular}{rrllllll} \toprule
                 & \multicolumn{6}{c}{Overlap}                            \\
      ranks & dofs    & $\sfrac{\delta}{H}=5\,\%$ & $t^{\text{walltime}}_{5\%}$ & $\sfrac{\delta}{H}=10\,\%$ & $t^{\text{walltime}}_{10\%}$ & $\sfrac{\delta}{H}=20\,\%$ & $t^{\text{walltime}}_{20\%}$\\ \midrule
      2   & $112\,964$    & $17$    & $ 61.4\;\unit{s}$ & $12$  & $ 63.1\;\unit{s}$ & $11$  & $ 65.5\;\unit{s}$ \\
      16  & $856\,168$    & $178$   & $ 97.9\;\unit{s}$ & $66$  & $ 98.6\;\unit{s}$ & $58$  & $127.0\;\unit{s}$ \\
      128 & $6\,661\,712$ & $>1000$ & $-$               & $200$ & $416.0\;\unit{s}$ & $171$ & $568.0\;\unit{s}$ \\\bottomrule
    \end{tabular}
  \end{center}
  \caption{
    \label{tab:itertation-count-maxwell-ybeamsplitter}
    Overview of the wall times and the number of GMRES iterations required to solve the three-dimensional y-beam splitter 
    using the ORAS preconditioner dependent on the number of ranks and the overlap between the subdomains.
  }
\end{table}

Our findings in Table~\ref{tab:itertation-count-maxwell-ybeamsplitter} indicate that a $10\,\%$ overlap is also sufficient for more complex problems. 
A larger overlap is required to ensure convergence even for even more complex 
problems, e.g., problems with a high amount of internal reflection.

\section{Conclusions}
In this work, we considered RAS preconditioners and optimized RAS preconditioners for the iterative GMRES linear solution within linear PDEs and within Newton's methods for solving 
nonlinear, coupled PDE systems.
Both classical Lagrange-based finite elements and N\'ed\'elec finite elements work in our implementation. 
The main technical work was the design, implementation, and debugging of the optimized Schwarz preconditioners 
in {\frosch} as well as of the {\tpetra}-based interface between the solver 
package {\frosch} and the finite element library {\dealii}. 

To demonstrate the performance, as model problems,
two challenging systems were chosen.
First, nonstationary, nonlinear FSI problems were considered.
The efficient solution of FSI remains an imminent problem, and open-source implementations 
in reliable, sustainable software packages remain rare. 
Due to the nonlinearities, a Newton solver with
line search was utilized, while the GMRES solvers with RAS preconditioners were used for the arising linear equation 
systems. For verification of the correct physics, the typical characteristic quantities of the benchmark problems were considered, and good agreement with the published literature was observed. The analysis of the RAS preconditioners showed 
good performance in two and three-dimensional numerical examples. 
Second, indefinite Maxwell's equations were considered, 
which remain challenging in their efficient and robust numerical solution. Again, we observed a very good performance of 
the ORAS preconditioner. The overall conclusion of this work 
is to have a very easy-to-use software interface 
while still giving satisfactory solver performances 
for challenging PDE problems.
Thus, we have obtained
very promising findings for the novel sustainable 
interface between {\frosch} and {\dealii} that can be further applied to other linear and nonlinear PDEs and coupled PDE systems.

\section*{Acknowledgments}
  This work is funded by the Deutsche Forschungsgemeinschaft (DFG) under Germany’s Excellence
  Strategy within the Cluster of Excellence PhoenixD (EXC 2122, Project ID 390833453).
  Furthermore, we would like to thank Jan Philipp Thiele for many fruitful discussions on the {\tpetra} interface of {\dealii}, 
  and Sven Beuchler for helpful discussions about preconditioners.

\bibliography{lit}


\begin{thebibliography}{81}


\ifx \showCODEN    \undefined \def \showCODEN     #1{\unskip}     \fi
\ifx \showDOI      \undefined \def \showDOI       #1{#1}\fi
\ifx \showISBNx    \undefined \def \showISBNx     #1{\unskip}     \fi
\ifx \showISBNxiii \undefined \def \showISBNxiii  #1{\unskip}     \fi
\ifx \showISSN     \undefined \def \showISSN      #1{\unskip}     \fi
\ifx \showLCCN     \undefined \def \showLCCN      #1{\unskip}     \fi
\ifx \shownote     \undefined \def \shownote      #1{#1}          \fi
\ifx \showarticletitle \undefined \def \showarticletitle #1{#1}   \fi
\ifx \showURL      \undefined \def \showURL       {\relax}        \fi
\providecommand\bibfield[2]{#2}
\providecommand\bibinfo[2]{#2}
\providecommand\natexlab[1]{#1}
\providecommand\showeprint[2][]{arXiv:#2}

\bibitem[Amestoy et~al\mbox{.}(2019)]%
        {Software:MUMPS}
\bibfield{author}{\bibinfo{person}{P.R. Amestoy}, \bibinfo{person}{A. Buttari},
  \bibinfo{person}{J.-Y. L'Excellent}, {and} \bibinfo{person}{T. Mary}.}
  \bibinfo{year}{2019}\natexlab{}.
\newblock \showarticletitle{{Performance and Scalability of the Block Low-Rank
  Multifrontal Factorization on Multicore Architectures}}.
\newblock \bibinfo{journal}{\emph{ACM Trans. Math. Software}}
  \bibinfo{volume}{45} (\bibinfo{year}{2019}), \bibinfo{pages}{2:1--2:26}.
\newblock
Issue 1.


\bibitem[Arndt et~al\mbox{.}(2023)]%
        {Software:dealii:95}
\bibfield{author}{\bibinfo{person}{D. Arndt}, \bibinfo{person}{W. Bangerth},
  \bibinfo{person}{M. Bergbauer}, \bibinfo{person}{M. Feder},
  \bibinfo{person}{M. Fehling}, \bibinfo{person}{J. Heinz}, \bibinfo{person}{T.
  Heister}, \bibinfo{person}{L. Heltai}, \bibinfo{person}{M. Kronbichler},
  \bibinfo{person}{M. Maier}, \bibinfo{person}{P. Munch},
  \bibinfo{person}{J.-P. Pelteret}, \bibinfo{person}{B. Turcksin},
  \bibinfo{person}{D. Wells}, {and} \bibinfo{person}{S. Zampini}.}
  \bibinfo{year}{2023}\natexlab{}.
\newblock \showarticletitle{The \texttt{deal.II} Library, Version 9.5}.
\newblock \bibinfo{journal}{\emph{Journal of Numerical Mathematics}}
  \bibinfo{volume}{31}, \bibinfo{number}{3} (\bibinfo{year}{2023}),
  \bibinfo{pages}{231--246}.
\newblock
\urldef\tempurl%
\url{https://doi.org/10.1515/jnma-2023-0089}
\showDOI{\tempurl}


\bibitem[Arndt et~al\mbox{.}(2021)]%
        {deal2020}
\bibfield{author}{\bibinfo{person}{D. Arndt}, \bibinfo{person}{W. Bangerth},
  \bibinfo{person}{D. Davydov}, \bibinfo{person}{T. Heister},
  \bibinfo{person}{L. Heltai}, \bibinfo{person}{M. Kronbichler},
  \bibinfo{person}{M. Maier}, \bibinfo{person}{J.-P. Pelteret},
  \bibinfo{person}{B. Turcksin}, {and} \bibinfo{person}{D. Wells}.}
  \bibinfo{year}{2021}\natexlab{}.
\newblock \showarticletitle{The {deal.II} finite element library: Design,
  features, and insights}.
\newblock \bibinfo{journal}{\emph{Computers {\&} Mathematics with
  Applications}}  \bibinfo{volume}{81} (\bibinfo{year}{2021}),
  \bibinfo{pages}{407--422}.
\newblock
\showISSN{0898-1221}
\urldef\tempurl%
\url{https://doi.org/10.1016/j.camwa.2020.02.022}
\showDOI{\tempurl}


\bibitem[Barker and Cai(2010)]%
        {barker_two-level_2010}
\bibfield{author}{\bibinfo{person}{A.~T. Barker} {and} \bibinfo{person}{X.-C.
  Cai}.} \bibinfo{year}{2010}\natexlab{}.
\newblock \showarticletitle{Two-{Level} {Newton} and {Hybrid} {Schwarz}
  {Preconditioners} for {Fluid}-{Structure} {Interaction}}.
\newblock \bibinfo{journal}{\emph{SIAM Journal on Scientific Computing}}
  \bibinfo{volume}{32}, \bibinfo{number}{4} (\bibinfo{date}{Jan.}
  \bibinfo{year}{2010}), \bibinfo{pages}{2395--2417}.
\newblock
\showISSN{1064-8275}
\urldef\tempurl%
\url{https://doi.org/10.1137/090779425}
\showDOI{\tempurl}
\newblock
\shownote{Publisher: Society for Industrial and Applied Mathematics}.


\bibitem[Bavier et~al\mbox{.}(2012)]%
        {bavier_amesos2_2012}
\bibfield{author}{\bibinfo{person}{E. Bavier}, \bibinfo{person}{M. Hoemmen},
  \bibinfo{person}{S. Rajamanickam}, {and} \bibinfo{person}{H. Thornquist}.}
  \bibinfo{year}{2012}\natexlab{}.
\newblock \showarticletitle{Amesos2 and {Belos}: {Direct} and iterative solvers
  for large sparse linear systems}.
\newblock \bibinfo{journal}{\emph{Scientific Programming}}
  \bibinfo{volume}{20}, \bibinfo{number}{3} (\bibinfo{date}{Jan.}
  \bibinfo{year}{2012}), \bibinfo{pages}{241--255}.
\newblock
\showISSN{1058-9244}
\urldef\tempurl%
\url{https://doi.org/10.3233/SPR-2012-0352}
\showDOI{\tempurl}
\newblock
\shownote{Publisher: IOS Press}.


\bibitem[Bazilevs et~al\mbox{.}(2013)]%
        {BaKeTe13}
\bibfield{author}{\bibinfo{person}{Y. Bazilevs}, \bibinfo{person}{K. Takizawa},
  {and} \bibinfo{person}{T.E. Tezduyar}.} \bibinfo{year}{2013}\natexlab{}.
\newblock \bibinfo{booktitle}{\emph{Computational Fluid-Structure Interaction:
  Methods and Applications}}.
\newblock \bibinfo{publisher}{Wiley}.
\newblock


\bibitem[Beuchler et~al\mbox{.}(2022)]%
        {Art:Kinnewig:22:DD26}
\bibfield{author}{\bibinfo{person}{S. Beuchler}, \bibinfo{person}{S. Kinnewig},
  {and} \bibinfo{person}{T. Wick}.} \bibinfo{year}{2022}\natexlab{}.
\newblock \showarticletitle{Parallel Domain Decomposition Solvers for the Time
  Harmonic Maxwell Equations}. In \bibinfo{booktitle}{\emph{Domain
  Decomposition Methods in Science and Engineering XXVI}},
  \bibfield{editor}{\bibinfo{person}{Susanne~C. Brenner}, \bibinfo{person}{Eric
  Chung}, \bibinfo{person}{Axel Klawonn}, \bibinfo{person}{Felix Kwok},
  \bibinfo{person}{Jinchao Xu}, {and} \bibinfo{person}{Jun Zou}} (Eds.).
  \bibinfo{publisher}{Springer International Publishing},
  \bibinfo{address}{Cham}, \bibinfo{pages}{653--660}.
\newblock
\showISBNx{978-3-030-95025-5}
\urldef\tempurl%
\url{https://doi.org/10.1007/978-3-030-95025-5_71}
\showDOI{\tempurl}


\bibitem[Boman et~al\mbox{.}(2012)]%
        {boman_zoltan_2012}
\bibfield{author}{\bibinfo{person}{E.~G. Boman}, \bibinfo{person}{Ü.~V.
  Çatalyürek}, \bibinfo{person}{C. Chevalier}, {and} \bibinfo{person}{K.~D.
  Devine}.} \bibinfo{year}{2012}\natexlab{}.
\newblock \showarticletitle{The {Zoltan} and {Isorropia} parallel toolkits for
  combinatorial scientific computing: {Partitioning}, ordering and coloring}.
\newblock \bibinfo{journal}{\emph{Sci. Program.}} \bibinfo{volume}{20},
  \bibinfo{number}{2} (\bibinfo{date}{April} \bibinfo{year}{2012}),
  \bibinfo{pages}{129--150}.
\newblock
\urldef\tempurl%
\url{https://doi.org/10.1155/2012/713587}
\showDOI{\tempurl}


\bibitem[Bootland et~al\mbox{.}(2024)]%
        {Art:Dolean:24::D27}
\bibfield{author}{\bibinfo{person}{N. Bootland}, \bibinfo{person}{S. Borzooei},
  \bibinfo{person}{V. Dolean}, {and} \bibinfo{person}{P.-H. Tournier}.}
  \bibinfo{year}{2024}\natexlab{}.
\newblock \showarticletitle{Numerical Assessment of PML Transmission Conditions
  in a Domain Decomposition Method for the Helmholtz Equation}. In
  \bibinfo{booktitle}{\emph{Domain Decomposition Methods in Science and
  Engineering XXVII}}, \bibfield{editor}{\bibinfo{person}{Zden{\v{e}}k
  Dost{\'a}l}, \bibinfo{person}{Tom{\'a}{\v{s}} Kozubek}, \bibinfo{person}{Axel
  Klawonn}, \bibinfo{person}{Ulrich Langer}, \bibinfo{person}{Luca~F.
  Pavarino}, \bibinfo{person}{Jakub {\v{S}}{\'i}stek}, {and}
  \bibinfo{person}{Olof~B. Widlund}} (Eds.). \bibinfo{publisher}{Springer
  Nature Switzerland}, \bibinfo{address}{Cham}, \bibinfo{pages}{445--453}.
\newblock
\urldef\tempurl%
\url{https://doi.org/10.1007/978-3-031-50769-4_53}
\showDOI{\tempurl}


\bibitem[Borzooei et~al\mbox{.}(2024)]%
        {Art:Dolean:TransmissionMaxwell:24}
\bibfield{author}{\bibinfo{person}{S. Borzooei}, \bibinfo{person}{V. Dolean},
  \bibinfo{person}{P.-H. Tournier}, {and} \bibinfo{person}{C. Migliaccio}.}
  \bibinfo{year}{2024}\natexlab{}.
\newblock \showarticletitle{Solution of Time-Harmonic Maxwell's Equations
  by a Domain Decomposition Method Based on PML Transmission Conditions}. In
  \bibinfo{booktitle}{\emph{Scientific Computing in Electrical Engineering}},
  \bibfield{editor}{\bibinfo{person}{Martijn van Beurden},
  \bibinfo{person}{Neil~V. Budko}, \bibinfo{person}{Gabriela Ciuprina},
  \bibinfo{person}{Wil Schilders}, \bibinfo{person}{Harshit Bansal}, {and}
  \bibinfo{person}{Ruxandra Barbulescu}} (Eds.). \bibinfo{publisher}{Springer
  Nature Switzerland}, \bibinfo{address}{Cham}, \bibinfo{pages}{45--52}.
\newblock
\urldef\tempurl%
\url{https://doi.org/10.1007/978-3-031-54517-7_5}
\showDOI{\tempurl}


\bibitem[Bungartz et~al\mbox{.}(2010)]%
        {BuSc10}
\bibfield{author}{\bibinfo{person}{H.-J. Bungartz}, \bibinfo{person}{M. Mehl},
  {and} \bibinfo{person}{M. Sch\"afer}.} \bibinfo{year}{2010}\natexlab{}.
\newblock \bibinfo{booktitle}{\emph{Fluid-Structure Interaction II: Modelling,
  Simulation, Optimization}}.
\newblock \bibinfo{publisher}{Springer}.
\newblock


\bibitem[Burstedde et~al\mbox{.}(2011)]%
        {Software:p4est:11}
\bibfield{author}{\bibinfo{person}{C. Burstedde}, \bibinfo{person}{L.~C.
  Wilcox}, {and} \bibinfo{person}{O. Ghattas}.}
  \bibinfo{year}{2011}\natexlab{}.
\newblock \showarticletitle{{\texttt{p4est}}: Scalable Algorithms for Parallel
  Adaptive Mesh Refinement on Forests of Octrees}.
\newblock \bibinfo{journal}{\emph{SIAM Journal on Scientific Computing}}
  \bibinfo{volume}{33}, \bibinfo{number}{3} (\bibinfo{year}{2011}),
  \bibinfo{pages}{1103--1133}.
\newblock
\urldef\tempurl%
\url{https://doi.org/10.1137/100791634}
\showDOI{\tempurl}


\bibitem[Cai and Sarkis(1999)]%
        {cai_restricted_1999}
\bibfield{author}{\bibinfo{person}{X.-C. Cai} {and} \bibinfo{person}{M.
  Sarkis}.} \bibinfo{year}{1999}\natexlab{}.
\newblock \showarticletitle{A {Restricted} {Additive} {Schwarz}
  {Preconditioner} for {General} {Sparse} {Linear} {Systems}}.
\newblock \bibinfo{journal}{\emph{SIAM Journal on Scientific Computing}}
  \bibinfo{volume}{21}, \bibinfo{number}{2} (\bibinfo{date}{Jan.}
  \bibinfo{year}{1999}), \bibinfo{pages}{792--797}.
\newblock
\showISSN{1064-8275}
\urldef\tempurl%
\url{https://doi.org/10.1137/S106482759732678X}
\showDOI{\tempurl}
\newblock
\shownote{Publisher: Society for Industrial and Applied Mathematics}.


\bibitem[Cros(2003)]%
        {cros_preconditioner_2003}
\bibfield{author}{\bibinfo{person}{J.-M. Cros}.}
  \bibinfo{year}{2003}\natexlab{}.
\newblock \showarticletitle{A preconditioner for the {Schur} complement domain
  decomposition method}. In \bibinfo{booktitle}{\emph{Domain {Decomposition}
  {Methods} in {Science} and {Engineering}}},
  \bibfield{editor}{\bibinfo{person}{O.~Widlund I.~Herrera, D.~Keyes} {and}
  \bibinfo{person}{R.~Yates}} (Eds.). \bibinfo{publisher}{National Autonomous
  University of Mexico (UNAM), Mexico City, Mexico, ISBN 970-32-0859-2},
  \bibinfo{pages}{373--380}.
\newblock


\bibitem[Crosetto et~al\mbox{.}(2011)]%
        {CrDeFouQua11}
\bibfield{author}{\bibinfo{person}{P. Crosetto}, \bibinfo{person}{S. Deparis},
  \bibinfo{person}{G. Fourestey}, {and} \bibinfo{person}{A. Quarteroni}.}
  \bibinfo{year}{2011}\natexlab{}.
\newblock \showarticletitle{Parallel algorithms for fluid-structure interaction
  problems in haemodynamics}.
\newblock \bibinfo{journal}{\emph{SIAM J. Sci. Comp.}} \bibinfo{volume}{33},
  \bibinfo{number}{4} (\bibinfo{year}{2011}), \bibinfo{pages}{1598--1622}.
\newblock


\bibitem[Deparis et~al\mbox{.}(2016)]%
        {DeparisFortiGrandperrinQuarteroni:2016}
\bibfield{author}{\bibinfo{person}{S. Deparis}, \bibinfo{person}{D. Forti},
  \bibinfo{person}{G. Grandperrin}, {and} \bibinfo{person}{A. Quarteroni}.}
  \bibinfo{year}{2016}\natexlab{}.
\newblock \showarticletitle{FaCSI: A block parallel preconditioner for
  fluid-structure interaction in hemodynamics}.
\newblock \bibinfo{journal}{\emph{J. Comput. Phys.}}  \bibinfo{volume}{327}
  (\bibinfo{year}{2016}), \bibinfo{pages}{700--718}.
\newblock


\bibitem[Deparis et~al\mbox{.}(2015)]%
        {deparis_comparison_2015}
\bibfield{author}{\bibinfo{person}{S. Deparis}, \bibinfo{person}{D. Forti},
  \bibinfo{person}{A. Heinlein}, \bibinfo{person}{A. Klawonn},
  \bibinfo{person}{A. Quarteroni}, {and} \bibinfo{person}{O. Rheinbach}.}
  \bibinfo{year}{2015}\natexlab{}.
\newblock \showarticletitle{A {Comparison} of {Preconditioners} for the
  {Steklov}–{Poincaré} {Formulation} of the {Fluid}-{Structure} {Coupling}
  in {Hemodynamics}}.
\newblock \bibinfo{journal}{\emph{PAMM}} \bibinfo{volume}{15},
  \bibinfo{number}{1} (\bibinfo{year}{2015}), \bibinfo{pages}{93--94}.
\newblock
\newblock
\shownote{Publisher: Wiley Online Library}.


\bibitem[Dohrmann(2003)]%
        {dohrmann_preconditioner_2003}
\bibfield{author}{\bibinfo{person}{C.~R. Dohrmann}.}
  \bibinfo{year}{2003}\natexlab{}.
\newblock \showarticletitle{A preconditioner for substructuring based on
  constrained energy minimization}.
\newblock \bibinfo{journal}{\emph{SIAM Journal on Scientific Computing}}
  \bibinfo{volume}{25}, \bibinfo{number}{1} (\bibinfo{year}{2003}),
  \bibinfo{pages}{246--258}.
\newblock
\showISSN{1064-8275}
\urldef\tempurl%
\url{https://doi.org/10.1137/S1064827502412887}
\showDOI{\tempurl}


\bibitem[Dohrmann et~al\mbox{.}(2008a)]%
        {dohrmann_domain_2008}
\bibfield{author}{\bibinfo{person}{C.~R. Dohrmann}, \bibinfo{person}{A.
  Klawonn}, {and} \bibinfo{person}{O.~B. Widlund}.}
  \bibinfo{year}{2008}\natexlab{a}.
\newblock \showarticletitle{Domain decomposition for less regular subdomains:
  overlapping {Schwarz} in two dimensions}.
\newblock \bibinfo{journal}{\emph{SIAM J. Numer. Anal.}} \bibinfo{volume}{46},
  \bibinfo{number}{4} (\bibinfo{year}{2008}), \bibinfo{pages}{2153--2168}.
\newblock
\showISSN{0036-1429}
\urldef\tempurl%
\url{https://doi.org/10.1137/070685841}
\showDOI{\tempurl}


\bibitem[Dohrmann et~al\mbox{.}(2008b)]%
        {dohrmann_family_2008}
\bibfield{author}{\bibinfo{person}{C.~R. Dohrmann}, \bibinfo{person}{A.
  Klawonn}, {and} \bibinfo{person}{O.~B. Widlund}.}
  \bibinfo{year}{2008}\natexlab{b}.
\newblock \showarticletitle{A family of energy minimizing coarse spaces for
  overlapping {Schwarz} preconditioners}.
\newblock In \bibinfo{booktitle}{\emph{Domain decomposition methods in science
  and engineering {XVII}}}. \bibinfo{series}{Lect. {Notes} {Comput}. {Sci}.
  {Eng}.}, Vol.~\bibinfo{volume}{60}. \bibinfo{publisher}{Springer, Berlin},
  \bibinfo{pages}{247--254}.
\newblock
\urldef\tempurl%
\url{https://doi.org/10.1007/978-3-540-75199-1_28}
\showDOI{\tempurl}


\bibitem[Dohrmann and Widlund(2017)]%
        {dohrmann_design_2017}
\bibfield{author}{\bibinfo{person}{C.~R. Dohrmann} {and} \bibinfo{person}{O.~B.
  Widlund}.} \bibinfo{year}{2017}\natexlab{}.
\newblock \showarticletitle{On the {Design} of {Small} {Coarse} {Spaces} for
  {Domain} {Decomposition} {Algorithms}}.
\newblock \bibinfo{journal}{\emph{SIAM Journal on Scientific Computing}}
  \bibinfo{volume}{39}, \bibinfo{number}{4} (\bibinfo{date}{Jan.}
  \bibinfo{year}{2017}), \bibinfo{pages}{A1466--A1488}.
\newblock
\showISSN{1064-8275}
\urldef\tempurl%
\url{https://doi.org/10.1137/17M1114272}
\showDOI{\tempurl}
\newblock
\shownote{Publisher: Society for Industrial and Applied Mathematics}.


\bibitem[Dolean et~al\mbox{.}(2009)]%
        {Art:DoGaGe:09}
\bibfield{author}{\bibinfo{person}{V. Dolean}, \bibinfo{person}{M.~J. Gander},
  {and} \bibinfo{person}{L. Gerardo-Giorda}.} \bibinfo{year}{2009}\natexlab{}.
\newblock \showarticletitle{Optimized Schwarz Methods for Maxwell's Equations}.
\newblock \bibinfo{journal}{\emph{SIAM Journal on Scientific Computing}}
  \bibinfo{volume}{31}, \bibinfo{number}{3} (\bibinfo{year}{2009}),
  \bibinfo{pages}{2193--2213}.
\newblock
\urldef\tempurl%
\url{https://doi.org/10.1137/080728536}
\showDOI{\tempurl}


\bibitem[Dolean et~al\mbox{.}(2015)]%
        {dolean_introduction_2015}
\bibfield{author}{\bibinfo{person}{V. Dolean}, \bibinfo{person}{P. Jolivet},
  {and} \bibinfo{person}{F. Nataf}.} \bibinfo{year}{2015}\natexlab{}.
\newblock \bibinfo{booktitle}{\emph{An {Introduction} to {Domain}
  {Decomposition} {Methods}: {Algorithms}, {Theory}, and {Parallel}
  {Implementation}}}.
\newblock \bibinfo{publisher}{Society for Industrial and Applied Mathematics},
  \bibinfo{address}{Philadelphia, PA}.
\newblock
\showISBNx{978-1-61197-405-8 978-1-61197-406-5}
\urldef\tempurl%
\url{https://doi.org/10.1137/1.9781611974065}
\showDOI{\tempurl}


\bibitem[Don\'ea et~al\mbox{.}(1977)]%
        {DoFaGi77}
\bibfield{author}{\bibinfo{person}{J. Don\'ea}, \bibinfo{person}{P.
  Fasoli-Stella}, {and} \bibinfo{person}{S. Giuliani}.}
  \bibinfo{year}{1977}\natexlab{}.
\newblock \showarticletitle{{L}agrangian and {E}ulerian finite element
  techniques for transient fluid-structure interaction problems}. In
  \bibinfo{booktitle}{\emph{Trans. 4th Int. Conf. on Structural Mechanics in
  Reactor Technology}}. \bibinfo{pages}{Paper B1/2}.
\newblock


\bibitem[Donea et~al\mbox{.}(2004)]%
        {DoHuePoRo04}
\bibfield{author}{\bibinfo{person}{J. Donea}, \bibinfo{person}{A. Huerta},
  \bibinfo{person}{J.-Ph. Ponthot}, {and} \bibinfo{person}{A.
  Rodriguez-Ferran}.} \bibinfo{year}{2004}\natexlab{}.
\newblock \bibinfo{booktitle}{\emph{Arbitrary {L}agrangian-{E}ulerian
  {m}ethods}}.
\newblock \bibinfo{publisher}{John Wiley and Sons}, \bibinfo{pages}{1--25}.
\newblock


\bibitem[Dunne(2006)]%
        {Du06}
\bibfield{author}{\bibinfo{person}{T. Dunne}.} \bibinfo{year}{2006}\natexlab{}.
\newblock \showarticletitle{An {E}ulerian approach to fluid-structure
  interaction and goal-oriented mesh adaption}.
\newblock \bibinfo{journal}{\emph{Int. J. Numer. Methods in Fluids}}
  \bibinfo{volume}{51} (\bibinfo{year}{2006}), \bibinfo{pages}{1017--1039}.
\newblock


\bibitem[Efstathiou and Gander(2003)]%
        {efstathiou_why_2003}
\bibfield{author}{\bibinfo{person}{E. Efstathiou} {and} \bibinfo{person}{M.~J.
  Gander}.} \bibinfo{year}{2003}\natexlab{}.
\newblock \showarticletitle{Why restricted additive {Schwarz} converges faster
  than additive {Schwarz}}.
\newblock \bibinfo{journal}{\emph{BIT. Numerical Mathematics}}
  \bibinfo{volume}{43}, \bibinfo{number}{suppl.} (\bibinfo{year}{2003}),
  \bibinfo{pages}{945--959}.
\newblock
\showISSN{0006-3835}
\urldef\tempurl%
\url{https://doi.org/10.1023/B:BITN.0000014563.33622.1d}
\showDOI{\tempurl}


\bibitem[{El Bouajaji} et~al\mbox{.}(2015)]%
        {Art:Bouajaji:15}
\bibfield{author}{\bibinfo{person}{M. {El Bouajaji}}, \bibinfo{person}{B.
  Thierry}, \bibinfo{person}{X. Antoine}, {and} \bibinfo{person}{C. Geuzaine}.}
  \bibinfo{year}{2015}\natexlab{}.
\newblock \showarticletitle{A quasi-optimal domain decomposition algorithm for
  the time-harmonic Maxwell's equations}.
\newblock \bibinfo{journal}{\emph{J. Comput. Phys.}}  \bibinfo{volume}{294}
  (\bibinfo{year}{2015}), \bibinfo{pages}{38--57}.
\newblock
\showISSN{0021-9991}
\urldef\tempurl%
\url{https://doi.org/10.1016/j.jcp.2015.03.041}
\showDOI{\tempurl}


\bibitem[Ernst and Gander(2012)]%
        {Art:Gander:Helmholtz:2012}
\bibfield{author}{\bibinfo{person}{O.~G. Ernst} {and} \bibinfo{person}{M.~J.
  Gander}.} \bibinfo{year}{2012}\natexlab{}.
\newblock \showarticletitle{Why it is difficult to solve {H}elmholtz problems
  with classical iterative methods}.
\newblock In \bibinfo{booktitle}{\emph{Numerical analysis of multiscale
  problems}}. \bibinfo{series}{Lect. Notes Comput. Sci. Eng.},
  Vol.~\bibinfo{volume}{83}. \bibinfo{publisher}{Springer, Heidelberg},
  \bibinfo{pages}{325--363}.
\newblock
\showISBNx{978-3-642-22060-9; 978-3-642-22061-6}
\urldef\tempurl%
\url{https://doi.org/10.1007/978-3-642-22061-6\_10}
\showDOI{\tempurl}


\bibitem[Farhat et~al\mbox{.}(2001)]%
        {farhat_feti-dp_2001}
\bibfield{author}{\bibinfo{person}{C. Farhat}, \bibinfo{person}{M. Lesoinne},
  \bibinfo{person}{P. LeTallec}, \bibinfo{person}{K. Pierson}, {and}
  \bibinfo{person}{D. Rixen}.} \bibinfo{year}{2001}\natexlab{}.
\newblock \showarticletitle{{FETI}-{DP}: a dual–primal unified {FETI}
  method—part {I}: {A} faster alternative to the two-level {FETI} method}.
\newblock \bibinfo{journal}{\emph{Internat. J. Numer. Methods Engrg.}}
  \bibinfo{volume}{50}, \bibinfo{number}{7} (\bibinfo{year}{2001}),
  \bibinfo{pages}{1523--1544}.
\newblock
\showISSN{1097-0207}
\urldef\tempurl%
\url{https://doi.org/10.1002/nme.76}
\showDOI{\tempurl}


\bibitem[Farhat et~al\mbox{.}(2000)]%
        {farhat_scalable_2000}
\bibfield{author}{\bibinfo{person}{C. Farhat}, \bibinfo{person}{M. Lesoinne},
  {and} \bibinfo{person}{K. Pierson}.} \bibinfo{year}{2000}\natexlab{}.
\newblock \showarticletitle{A scalable dual-primal domain decomposition
  method}.
\newblock \bibinfo{journal}{\emph{Numerical Linear Algebra with Applications}}
  \bibinfo{volume}{7}, \bibinfo{number}{7-8} (\bibinfo{year}{2000}),
  \bibinfo{pages}{687--714}.
\newblock
\showISSN{1099-1506}
\urldef\tempurl%
\url{https://doi.org/10.1002/1099-1506(200010/12)7:7/8<687::AID-NLA219>3.0.CO;2-S}
\showDOI{\tempurl}


\bibitem[Feynman et~al\mbox{.}(1963)]%
        {Bk:Fey:63}
\bibfield{author}{\bibinfo{person}{R.~P. Feynman}, \bibinfo{person}{R.~B.
  Leighton}, {and} \bibinfo{person}{M. Sands}.}
  \bibinfo{year}{1963}\natexlab{}.
\newblock \bibinfo{booktitle}{\emph{The Feynman Lectures on Physics. Vol. II.
  Mainly Electromagnetism and Matter}}.
\newblock \bibinfo{publisher}{California Institute of Technology, Michael A.
  Gottlieb and Rudolf Pfeiffer}.
\newblock


\bibitem[Formaggia and Nobile(1999)]%
        {FoNo99}
\bibfield{author}{\bibinfo{person}{L. Formaggia} {and} \bibinfo{person}{F.
  Nobile}.} \bibinfo{year}{1999}\natexlab{}.
\newblock \showarticletitle{A stability analysis for the Arbitrary {L}agrangian
  {E}ulerian Formulation with Finite Elements}.
\newblock \bibinfo{journal}{\emph{East-West Journal of Numerical Mathematics}}
  \bibinfo{volume}{7} (\bibinfo{year}{1999}), \bibinfo{pages}{105 -- 132}.
\newblock


\bibitem[Forti et~al\mbox{.}(2017)]%
        {FortiQuarteroniDeparis:2017}
\bibfield{author}{\bibinfo{person}{D. Forti}, \bibinfo{person}{A. Quarteroni},
  {and} \bibinfo{person}{S. Deparis}.} \bibinfo{year}{2017}\natexlab{}.
\newblock \showarticletitle{A parallel algorithm for the solution of
  large-scale nonconforming fluid-structure interaction problems in
  hemodynamics}.
\newblock \bibinfo{journal}{\emph{J. Comput. Math.}} \bibinfo{volume}{35},
  \bibinfo{number}{3} (\bibinfo{year}{2017}), \bibinfo{pages}{363--380}.
\newblock
\showISSN{0254-9409}


\bibitem[Galdi and Rannacher(2010)]%
        {GaRa10}
\bibfield{author}{\bibinfo{person}{G. Galdi} {and} \bibinfo{person}{R.
  Rannacher}.} \bibinfo{year}{2010}\natexlab{}.
\newblock \bibinfo{booktitle}{\emph{Fundamental Trends in Fluid-Structure
  Interaction}}.
\newblock \bibinfo{publisher}{World Scientific}. 293 pages.
\newblock


\bibitem[Gander(2006)]%
        {Art:Gander:06}
\bibfield{author}{\bibinfo{person}{M.~J. Gander}.}
  \bibinfo{year}{2006}\natexlab{}.
\newblock \showarticletitle{Optimized Schwarz Methods}.
\newblock \bibinfo{journal}{\emph{SIAM J. Numer. Anal.}} \bibinfo{volume}{44},
  \bibinfo{number}{2} (\bibinfo{year}{2006}), \bibinfo{pages}{699--731}.
\newblock
\urldef\tempurl%
\url{https://doi.org/10.1137/S0036142903425409}
\showDOI{\tempurl}


\bibitem[Gander and Zhang(2016)]%
        {Art:Gander:16}
\bibfield{author}{\bibinfo{person}{M.~J. Gander} {and} \bibinfo{person}{H.
  Zhang}.} \bibinfo{year}{2016}\natexlab{}.
\newblock \showarticletitle{Optimized Schwarz Methods with Overlap for the
  Helmholtz Equation}.
\newblock \bibinfo{journal}{\emph{SIAM Journal on Scientific Computing}}
  \bibinfo{volume}{38}, \bibinfo{number}{5} (\bibinfo{year}{2016}),
  \bibinfo{pages}{A3195--A3219}.
\newblock
\urldef\tempurl%
\url{https://doi.org/10.1137/15M1021659}
\showDOI{\tempurl}


\bibitem[Gee et~al\mbox{.}(2011)]%
        {gee2011truly}
\bibfield{author}{\bibinfo{person}{M.~W. Gee}, \bibinfo{person}{U.
  K{\"u}ttler}, {and} \bibinfo{person}{W.~A. Wall}.}
  \bibinfo{year}{2011}\natexlab{}.
\newblock \showarticletitle{Truly monolithic algebraic multigrid for
  fluid--structure interaction}.
\newblock \bibinfo{journal}{\emph{Int. J. Numer. Meth. Engrg.}}
  \bibinfo{volume}{85}, \bibinfo{number}{8} (\bibinfo{year}{2011}),
  \bibinfo{pages}{987--1016}.
\newblock


\bibitem[Goll et~al\mbox{.}(2017)]%
        {DOpElib}
\bibfield{author}{\bibinfo{person}{C. Goll}, \bibinfo{person}{T. Wick}, {and}
  \bibinfo{person}{W. Wollner}.} \bibinfo{year}{2017}\natexlab{}.
\newblock \showarticletitle{{DOpElib: Differential Equations and Optimization
  Environment; A Goal Oriented Software Library for Solving PDEs and
  Optimization Problems with PDEs}}.
\newblock \bibinfo{journal}{\emph{Archive of Numerical Software}}
  \bibinfo{volume}{5}, \bibinfo{number}{2} (\bibinfo{year}{2017}),
  \bibinfo{pages}{1--14}.
\newblock
\urldef\tempurl%
\url{https://doi.org/10.11588/ans.2017.2.11815}
\showDOI{\tempurl}


\bibitem[Grayver and Kolev(2015)]%
        {Art:Grayver:Geoelectromagnetic:2015}
\bibfield{author}{\bibinfo{person}{A.~V. Grayver} {and} \bibinfo{person}{T.~V.
  Kolev}.} \bibinfo{year}{2015}\natexlab{}.
\newblock \showarticletitle{Large-Scale {{3D}} Geoelectromagnetic Modeling
  Using Parallel Adaptive High-Order Finite Element Method}.
\newblock \bibinfo{journal}{\emph{Geophysics}} \bibinfo{volume}{80},
  \bibinfo{number}{6} (\bibinfo{year}{2015}), \bibinfo{pages}{E277--E291}.
\newblock
\showISSN{0016-8033, 1942-2156}
\urldef\tempurl%
\url{https://doi.org/10.1190/geo2015-0013.1}
\showDOI{\tempurl}


\bibitem[Heil(2004)]%
        {He04}
\bibfield{author}{\bibinfo{person}{M. Heil}.} \bibinfo{year}{2004}\natexlab{}.
\newblock \showarticletitle{An efficient solver for the fully coupled solution
  of large-displacement fluid-structure interaction problems}.
\newblock \bibinfo{journal}{\emph{Comput. Methods Appl. Mech. Engrg.}}
  \bibinfo{volume}{193} (\bibinfo{year}{2004}), \bibinfo{pages}{1--23}.
\newblock


\bibitem[Heinlein(2016)]%
        {heinlein_parallel_2016-2}
\bibfield{author}{\bibinfo{person}{A. Heinlein}.}
  \bibinfo{year}{2016}\natexlab{}.
\newblock \emph{\bibinfo{title}{Parallel overlapping {Schwarz} preconditioners
  and multiscale discretizations with applications to fluid-structure
  interaction and highly heterogeneous problems}}.
\newblock {PhD} {Thesis}. \bibinfo{school}{Universität zu Köln}.
\newblock


\bibitem[Heinlein et~al\mbox{.}(2019a)]%
        {heinlein_monolithic_2019}
\bibfield{author}{\bibinfo{person}{A. Heinlein}, \bibinfo{person}{C. Hochmuth},
  {and} \bibinfo{person}{A. Klawonn}.} \bibinfo{year}{2019}\natexlab{a}.
\newblock \showarticletitle{Monolithic overlapping {Schwarz} domain
  decomposition methods with {GDSW} coarse spaces for incompressible fluid flow
  problems}.
\newblock \bibinfo{journal}{\emph{SIAM Journal on Scientific Computing}}
  \bibinfo{volume}{41}, \bibinfo{number}{4} (\bibinfo{year}{2019}),
  \bibinfo{pages}{C291--C316}.
\newblock
\showISSN{1064-8275}
\urldef\tempurl%
\url{https://doi.org/10.1137/18M1184047}
\showDOI{\tempurl}


\bibitem[Heinlein et~al\mbox{.}(2021)]%
        {heinlein_fully_2021}
\bibfield{author}{\bibinfo{person}{A. Heinlein}, \bibinfo{person}{C. Hochmuth},
  {and} \bibinfo{person}{A. Klawonn}.} \bibinfo{year}{2021}\natexlab{}.
\newblock \showarticletitle{Fully algebraic two-level overlapping {Schwarz}
  preconditioners for elasticity problems}. In
  \bibinfo{booktitle}{\emph{Numerical mathematics and advanced
  applications—{ENUMATH} 2019}} \emph{(\bibinfo{series}{Lect. {Notes}
  {Comput}. {Sci}. {Eng}.}, Vol.~\bibinfo{volume}{139})}.
  \bibinfo{publisher}{Springer, Cham}, \bibinfo{pages}{531--539}.
\newblock
\urldef\tempurl%
\url{https://doi.org/10.1007/978-3-030-55874-1_52}
\showDOI{\tempurl}


\bibitem[Heinlein et~al\mbox{.}(2019b)]%
        {heinlein_adaptive_2019}
\bibfield{author}{\bibinfo{person}{A. Heinlein}, \bibinfo{person}{A. Klawonn},
  \bibinfo{person}{J. Knepper}, {and} \bibinfo{person}{O. Rheinbach}.}
  \bibinfo{year}{2019}\natexlab{b}.
\newblock \showarticletitle{Adaptive {GDSW} coarse spaces for overlapping
  {Schwarz} methods in three dimensions}.
\newblock \bibinfo{journal}{\emph{SIAM Journal on Scientific Computing}}
  \bibinfo{volume}{41}, \bibinfo{number}{5} (\bibinfo{year}{2019}),
  \bibinfo{pages}{A3045--A3072}.
\newblock
\showISSN{1064-8275}
\urldef\tempurl%
\url{https://doi.org/10.1137/18M1220613}
\showDOI{\tempurl}


\bibitem[Heinlein et~al\mbox{.}(2020)]%
        {heinlein_frosch_2020}
\bibfield{author}{\bibinfo{person}{A. Heinlein}, \bibinfo{person}{A. Klawonn},
  \bibinfo{person}{S. Rajamanickam}, {and} \bibinfo{person}{O. Rheinbach}.}
  \bibinfo{year}{2020}\natexlab{}.
\newblock \showarticletitle{{FROSch}: a fast and robust overlapping {Schwarz}
  domain decomposition preconditioner based on {Xpetra} in {Trilinos}}. In
  \bibinfo{booktitle}{\emph{Domain decomposition methods in science and
  engineering {XXV}}} \emph{(\bibinfo{series}{Lect. {Notes} {Comput}. {Sci}.
  {Eng}.}, Vol.~\bibinfo{volume}{138})}. \bibinfo{publisher}{Springer, Cham},
  \bibinfo{pages}{176--184}.
\newblock
\urldef\tempurl%
\url{https://doi.org/10.1007/978-3-030-56750-7_19}
\showDOI{\tempurl}


\bibitem[Heinlein et~al\mbox{.}(2016a)]%
        {heinlein_parallel_2016-1}
\bibfield{author}{\bibinfo{person}{A. Heinlein}, \bibinfo{person}{A. Klawonn},
  {and} \bibinfo{person}{O. Rheinbach}.} \bibinfo{year}{2016}\natexlab{a}.
\newblock \showarticletitle{A parallel implementation of a two-level
  overlapping {Schwarz} method with energy-minimizing coarse space based on
  {Trilinos}}.
\newblock \bibinfo{journal}{\emph{SIAM Journal on Scientific Computing}}
  \bibinfo{volume}{38}, \bibinfo{number}{6} (\bibinfo{year}{2016}),
  \bibinfo{pages}{C713--C747}.
\newblock
\showISSN{1064-8275}
\urldef\tempurl%
\url{https://doi.org/10.1137/16M1062843}
\showDOI{\tempurl}


\bibitem[Heinlein et~al\mbox{.}(2016b)]%
        {heinlein_parallel_2016}
\bibfield{author}{\bibinfo{person}{A. Heinlein}, \bibinfo{person}{A. Klawonn},
  {and} \bibinfo{person}{O. Rheinbach}.} \bibinfo{year}{2016}\natexlab{b}.
\newblock \showarticletitle{Parallel two-level overlapping {Schwarz} methods in
  fluid-structure interaction}. In \bibinfo{booktitle}{\emph{Numerical
  mathematics and advanced applications—{ENUMATH} 2015}}
  \emph{(\bibinfo{series}{Lect. {Notes} {Comput}. {Sci}. {Eng}.},
  Vol.~\bibinfo{volume}{112})}. \bibinfo{publisher}{Springer, [Cham]},
  \bibinfo{pages}{521--530}.
\newblock
\urldef\tempurl%
\url{https://doi.org/10.1007/978-3-319-39929-4_5}
\showDOI{\tempurl}


\bibitem[Heinlein et~al\mbox{.}(2018)]%
        {heinlein_improving_2018}
\bibfield{author}{\bibinfo{person}{A. Heinlein}, \bibinfo{person}{A. Klawonn},
  \bibinfo{person}{O. Rheinbach}, {and} \bibinfo{person}{O.~B. Widlund}.}
  \bibinfo{year}{2018}\natexlab{}.
\newblock \showarticletitle{Improving the parallel performance of overlapping
  {Schwarz} methods by using a smaller energy minimizing coarse space}. In
  \bibinfo{booktitle}{\emph{Domain decomposition methods in science and
  engineering {XXIV}}} \emph{(\bibinfo{series}{Lect. {Notes} {Comput}. {Sci}.
  {Eng}.}, Vol.~\bibinfo{volume}{125})}. \bibinfo{publisher}{Springer, Cham},
  \bibinfo{pages}{383--392}.
\newblock
\urldef\tempurl%
\url{https://doi.org/10.1007/978-3-319-93873-8_3}
\showDOI{\tempurl}


\bibitem[Heinlein et~al\mbox{.}(2022a)]%
        {heinlein_frosch_2022}
\bibfield{author}{\bibinfo{person}{A. Heinlein}, \bibinfo{person}{M. Perego},
  {and} \bibinfo{person}{S. Rajamanickam}.} \bibinfo{year}{2022}\natexlab{a}.
\newblock \showarticletitle{{FROSch} {Preconditioners} for {Land} {Ice}
  {Simulations} of {Greenland} and {Antarctica}}.
\newblock \bibinfo{journal}{\emph{SIAM Journal on Scientific Computing}}
  \bibinfo{volume}{44}, \bibinfo{number}{2} (\bibinfo{date}{April}
  \bibinfo{year}{2022}), \bibinfo{pages}{B339--B367}.
\newblock
\showISSN{1064-8275}
\urldef\tempurl%
\url{https://doi.org/10.1137/21M1395260}
\showDOI{\tempurl}
\newblock
\shownote{Publisher: Society for Industrial and Applied Mathematics}.


\bibitem[Heinlein et~al\mbox{.}(2022b)]%
        {heinlein_parallel_2022}
\bibfield{author}{\bibinfo{person}{A. Heinlein}, \bibinfo{person}{O.
  Rheinbach}, {and} \bibinfo{person}{F. Röver}.}
  \bibinfo{year}{2022}\natexlab{b}.
\newblock \showarticletitle{Parallel {Scalability} of {Three}-{Level} {FROSch}
  {Preconditioners} to 220000 {Cores} using the {Theta} {Supercomputer}}.
\newblock \bibinfo{journal}{\emph{SIAM Journal on Scientific Computing}}
  (\bibinfo{date}{Aug.} \bibinfo{year}{2022}), \bibinfo{pages}{S173--S198}.
\newblock
\showISSN{1064-8275}
\urldef\tempurl%
\url{https://doi.org/10.1137/21M1431205}
\showDOI{\tempurl}
\newblock
\shownote{Publisher: Society for Industrial and Applied Mathematics}.


\bibitem[Heinlein et~al\mbox{.}(2023)]%
        {heinlein_multilevel_2023}
\bibfield{author}{\bibinfo{person}{A. Heinlein}, \bibinfo{person}{O.
  Rheinbach}, {and} \bibinfo{person}{F. Röver}.}
  \bibinfo{year}{2023}\natexlab{}.
\newblock \showarticletitle{A {Multilevel} {Extension} of the {GDSW}
  {Overlapping} {Schwarz} {Preconditioner} in {Two} {Dimensions}}.
\newblock \bibinfo{journal}{\emph{Computational Methods in Applied
  Mathematics}} \bibinfo{volume}{23}, \bibinfo{number}{4} (\bibinfo{date}{Oct.}
  \bibinfo{year}{2023}), \bibinfo{pages}{953--968}.
\newblock
\showISSN{1609-9389}
\urldef\tempurl%
\url{https://doi.org/10.1515/cmam-2022-0168}
\showDOI{\tempurl}
\newblock
\shownote{Publisher: De Gruyter}.


\bibitem[Heister and Wick(2020)]%
        {HeiWi20}
\bibfield{author}{\bibinfo{person}{T. Heister} {and} \bibinfo{person}{T.
  Wick}.} \bibinfo{year}{2020}\natexlab{}.
\newblock \showarticletitle{pfm-cracks: A parallel-adaptive framework for
  phase-field fracture propagation}.
\newblock \bibinfo{journal}{\emph{Software Impacts}}  \bibinfo{volume}{6}
  (\bibinfo{year}{2020}), \bibinfo{pages}{100045}.
\newblock
\showISSN{2665-9638}
\urldef\tempurl%
\url{https://doi.org/10.1016/j.simpa.2020.100045}
\showDOI{\tempurl}


\bibitem[Heroux et~al\mbox{.}(2005)]%
        {heroux_overview_2005}
\bibfield{author}{\bibinfo{person}{M.~A. Heroux}, \bibinfo{person}{R.~A.
  Bartlett}, \bibinfo{person}{V.~E. Howle}, \bibinfo{person}{R.~J. Hoekstra},
  \bibinfo{person}{J.~J. Hu}, \bibinfo{person}{T.~G. Kolda},
  \bibinfo{person}{R.~B. Lehoucq}, \bibinfo{person}{K.~R. Long},
  \bibinfo{person}{R.~P. Pawlowski}, \bibinfo{person}{E.~T. Phipps},
  \bibinfo{person}{A.~G. Salinger}, \bibinfo{person}{H.~K. Thornquist},
  \bibinfo{person}{R.~S. Tuminaro}, \bibinfo{person}{J.~M. Willenbring},
  \bibinfo{person}{A. Williams}, {and} \bibinfo{person}{K.~S. Stanley}.}
  \bibinfo{year}{2005}\natexlab{}.
\newblock \showarticletitle{An overview of the {Trilinos} {Project}}.
\newblock \bibinfo{journal}{\emph{Association for Computing Machinery.
  Transactions on Mathematical Software}} \bibinfo{volume}{31},
  \bibinfo{number}{3} (\bibinfo{year}{2005}), \bibinfo{pages}{397--423}.
\newblock
\showISSN{0098-3500}
\urldef\tempurl%
\url{https://doi.org/10.1145/1089014.1089021}
\showDOI{\tempurl}


\bibitem[Heywood et~al\mbox{.}(1996)]%
        {HeRaTu96}
\bibfield{author}{\bibinfo{person}{J.~G. Heywood}, \bibinfo{person}{R.
  Rannacher}, {and} \bibinfo{person}{S. Turek}.}
  \bibinfo{year}{1996}\natexlab{}.
\newblock \showarticletitle{Artificial boundaries and flux and pressure
  conditions for the incompressible {N}avier-{S}tokes Equations}.
\newblock \bibinfo{journal}{\emph{International Journal of Numerical Methods in
  Fluids}}  \bibinfo{volume}{22} (\bibinfo{year}{1996}),
  \bibinfo{pages}{325--352}.
\newblock


\bibitem[Hron and Turek(2006)]%
        {HrTu06b}
\bibfield{author}{\bibinfo{person}{J. Hron} {and} \bibinfo{person}{S. Turek}.}
  \bibinfo{year}{2006}\natexlab{}.
\newblock \bibinfo{booktitle}{\emph{Proposal for numerical benchmarking of
  fluid-structure interaction between an elastic object and laminar
  incompressible flow}}. Vol.~\bibinfo{volume}{53}.
\newblock \bibinfo{publisher}{Springer-Verlag}, \bibinfo{pages}{146 -- 170}.
\newblock


\bibitem[Hughes et~al\mbox{.}(1981)]%
        {HuLiZi81}
\bibfield{author}{\bibinfo{person}{T.J.R. Hughes}, \bibinfo{person}{W.K. Liu},
  {and} \bibinfo{person}{T Zimmermann}.} \bibinfo{year}{1981}\natexlab{}.
\newblock \showarticletitle{{L}agrangian-{E}ulerian finite element formulation
  for incompressible viscous flows}.
\newblock \bibinfo{journal}{\emph{Comput. Methods Appl. Mech. Engrg.}}
  \bibinfo{volume}{29} (\bibinfo{year}{1981}), \bibinfo{pages}{329--349}.
\newblock


\bibitem[Huynh et~al\mbox{.}(2024)]%
        {huynh_gdsw_2024}
\bibfield{author}{\bibinfo{person}{N.~M.~M. Huynh}, \bibinfo{person}{L.~F.
  Pavarino}, {and} \bibinfo{person}{S. Scacchi}.}
  \bibinfo{year}{2024}\natexlab{}.
\newblock \bibinfo{title}{{GDSW} preconditioners for composite {Discontinuous}
  {Galerkin} discretizations of multicompartment reaction-diffusion problems}.
\newblock
\newblock
\urldef\tempurl%
\url{https://doi.org/10.48550/arXiv.2405.17601}
\showDOI{\tempurl}
\newblock
\shownote{arXiv:2405.17601 [cs, math]}.


\bibitem[Jodlbauer et~al\mbox{.}(2019)]%
        {JoLaWi19_fsi}
\bibfield{author}{\bibinfo{person}{D. Jodlbauer}, \bibinfo{person}{U. Langer},
  {and} \bibinfo{person}{T. Wick}.} \bibinfo{year}{2019}\natexlab{}.
\newblock \showarticletitle{Parallel block-preconditioned monolithic solvers
  for fluid-structure interaction problems}.
\newblock \bibinfo{journal}{\emph{Int. J. Num. Meth. Eng.}}
  \bibinfo{volume}{117}, \bibinfo{number}{6} (\bibinfo{year}{2019}),
  \bibinfo{pages}{623--643}.
\newblock


\bibitem[Karypis and Kumar(1998)]%
        {karypis_fast_1998}
\bibfield{author}{\bibinfo{person}{G. Karypis} {and} \bibinfo{person}{V.
  Kumar}.} \bibinfo{year}{1998}\natexlab{}.
\newblock \showarticletitle{A {Fast} and {High} {Quality} {Multilevel} {Scheme}
  for {Partitioning} {Irregular} {Graphs}}.
\newblock \bibinfo{journal}{\emph{SIAM Journal on Scientific Computing}}
  \bibinfo{volume}{20}, \bibinfo{number}{1} (\bibinfo{date}{Jan.}
  \bibinfo{year}{1998}), \bibinfo{pages}{359--392}.
\newblock
\showISSN{1064-8275}
\urldef\tempurl%
\url{https://doi.org/10.1137/S1064827595287997}
\showDOI{\tempurl}
\newblock
\shownote{Publisher: Society for Industrial and Applied Mathematics}.


\bibitem[Kinnewig et~al\mbox{.}(2023)]%
        {Art:Kinnewig:23}
\bibfield{author}{\bibinfo{person}{S. Kinnewig}, \bibinfo{person}{T. Wick},
  {and} \bibinfo{person}{S. Beuchler}.} \bibinfo{year}{2023}\natexlab{}.
\newblock \showarticletitle{Algorithmic realization of the solution to the sign
  conflict problem for hanging nodes on hp-hexahedral N\'ed\'elec elements}.
\newblock  (\bibinfo{year}{2023}).
\newblock
\showeprint[arxiv]{2306.01416}~[math.NA]


\bibitem[Lions(1988)]%
        {lions_schwarz_1988}
\bibfield{author}{\bibinfo{person}{P.-L. Lions}.}
  \bibinfo{year}{1988}\natexlab{}.
\newblock \showarticletitle{On the {Schwarz} alternating method. {I}}.
\newblock In \bibinfo{booktitle}{\emph{First {International} {Symposium} on
  {Domain} {Decomposition} {Methods} for {Partial} {Differential} {Equations}
  ({Paris}, 1987)}}. \bibinfo{publisher}{SIAM, Philadelphia, PA},
  \bibinfo{pages}{1--42}.
\newblock
\urldef\tempurl%
\url{https://mathscinet.ams.org/mathscinet-getitem?mr=972510}
\showURL{%
\tempurl}


\bibitem[Mahmudlu et~al\mbox{.}(2023)]%
        {Art:Kues:23}
\bibfield{author}{\bibinfo{person}{H. Mahmudlu}, \bibinfo{person}{R.
  Johanning}, \bibinfo{person}{A. Van~Rees}, \bibinfo{person}{A.
  Khodadad~Kashi}, \bibinfo{person}{J.~P. Epping}, \bibinfo{person}{R. Haldar},
  \bibinfo{person}{K.-J. Boller}, {and} \bibinfo{person}{M. Kues}.}
  \bibinfo{year}{2023}\natexlab{}.
\newblock \showarticletitle{Fully On-Chip Photonic Turnkey Quantum Source for
  Entangled Qubit/Qudit State Generation}.
\newblock \bibinfo{journal}{\emph{Nature Photonics}} (\bibinfo{year}{2023}).
\newblock
\showISSN{1749-4885, 1749-4893}
\urldef\tempurl%
\url{https://doi.org/10.1038/s41566-023-01193-1}
\showDOI{\tempurl}


\bibitem[Melchert et~al\mbox{.}(2023)]%
        {Art:Melchert:Soliton:2023}
\bibfield{author}{\bibinfo{person}{O. Melchert}, \bibinfo{person}{S. Kinnewig},
  \bibinfo{person}{F. Dencker}, \bibinfo{person}{D. Perevoznik},
  \bibinfo{person}{S. Willms}, \bibinfo{person}{I.~V. Babushkin},
  \bibinfo{person}{M.~C. Wurz}, \bibinfo{person}{M. Kues}, \bibinfo{person}{S.
  Beuchler}, \bibinfo{person}{T. Wick}, \bibinfo{person}{U. Morgner}, {and}
  \bibinfo{person}{A. Demircan}.} \bibinfo{year}{2023}\natexlab{}.
\newblock \showarticletitle{Soliton Compression and Supercontinuum Spectra in
  Nonlinear Diamond Photonics}.
\newblock \bibinfo{journal}{\emph{Diamond and Related Materials}}
  \bibinfo{volume}{136} (\bibinfo{year}{2023}), \bibinfo{pages}{109939}.
\newblock
\showISSN{0925-9635}
\urldef\tempurl%
\url{https://doi.org/10.1016/j.diamond.2023.109939}
\showDOI{\tempurl}


\bibitem[Monk(2003)]%
        {Bk:Monk:2003}
\bibfield{author}{\bibinfo{person}{P. Monk}.} \bibinfo{year}{2003}\natexlab{}.
\newblock \bibinfo{booktitle}{\emph{Finite Element Methods for {{Maxwell}}'s
  Equations}}.
\newblock \bibinfo{publisher}{{Clarendon Press ; Oxford University Press}},
  \bibinfo{address}{{Oxford : New York}}.
\newblock
\showISBNx{978-0-19-850888-5}
\showLCCN{QC760 .M56 2003}


\bibitem[Nicolaides(1987)]%
        {nicolaides_deflation_1987}
\bibfield{author}{\bibinfo{person}{R.~A. Nicolaides}.}
  \bibinfo{year}{1987}\natexlab{}.
\newblock \showarticletitle{Deflation of conjugate gradients with applications
  to boundary value problems}.
\newblock \bibinfo{journal}{\emph{SIAM J. Numer. Anal.}} \bibinfo{volume}{24},
  \bibinfo{number}{2} (\bibinfo{year}{1987}), \bibinfo{pages}{355--365}.
\newblock
\showISSN{0036-1429}
\urldef\tempurl%
\url{https://doi.org/10.1137/0724027}
\showDOI{\tempurl}


\bibitem[Quarteroni and Valli(1999)]%
        {quarteroni_domain_1999}
\bibfield{author}{\bibinfo{person}{A. Quarteroni} {and} \bibinfo{person}{A.
  Valli}.} \bibinfo{year}{1999}\natexlab{}.
\newblock \bibinfo{booktitle}{\emph{Domain {Decomposition} {Methods} for
  {Partial} {Differential} {Equations}}}.
\newblock \bibinfo{publisher}{Oxford University Press},
  \bibinfo{address}{Oxford, New York}.
\newblock
\showISBNx{978-0-19-850178-7}


\bibitem[Rajamanickam et~al\mbox{.}(2021)]%
        {rajamanickam_kokkos_2021}
\bibfield{author}{\bibinfo{person}{S. Rajamanickam}, \bibinfo{person}{S. Acer},
  \bibinfo{person}{L. Berger-Vergiat}, \bibinfo{person}{V. Dang},
  \bibinfo{person}{N. Ellingwood}, \bibinfo{person}{E. Harvey},
  \bibinfo{person}{B. Kelley}, \bibinfo{person}{C.~R. Trott},
  \bibinfo{person}{J. Wilke}, {and} \bibinfo{person}{I. Yamazaki}.}
  \bibinfo{year}{2021}\natexlab{}.
\newblock \bibinfo{title}{Kokkos {Kernels}: {Performance} {Portable}
  {Sparse}/{Dense} {Linear} {Algebra} and {Graph} {Kernels}}.
\newblock
\newblock
\urldef\tempurl%
\url{https://doi.org/10.48550/arXiv.2103.11991}
\showDOI{\tempurl}
\newblock
\shownote{arXiv:2103.11991 [cs]}.


\bibitem[Richter(2017)]%
        {Ri17_fsi}
\bibfield{author}{\bibinfo{person}{T. Richter}.}
  \bibinfo{year}{2017}\natexlab{}.
\newblock \bibinfo{booktitle}{\emph{Fluid-structure interactions: models,
  analysis, and finite elements}}.
\newblock \bibinfo{publisher}{Springer}.
\newblock


\bibitem[Saad and Schultz(1986)]%
        {saad_gmres_1986}
\bibfield{author}{\bibinfo{person}{Y. Saad} {and} \bibinfo{person}{M.~H.
  Schultz}.} \bibinfo{year}{1986}\natexlab{}.
\newblock \showarticletitle{{GMRES}: a generalized minimal residual algorithm
  for solving nonsymmetric linear systems}.
\newblock \bibinfo{journal}{\emph{Society for Industrial and Applied
  Mathematics. Journal on Scientific and Statistical Computing}}
  \bibinfo{volume}{7}, \bibinfo{number}{3} (\bibinfo{year}{1986}),
  \bibinfo{pages}{856--869}.
\newblock
\showISSN{0196-5204}
\urldef\tempurl%
\url{https://doi.org/10.1137/0907058}
\showDOI{\tempurl}


\bibitem[Smith et~al\mbox{.}(1996)]%
        {smith_domain_1996}
\bibfield{author}{\bibinfo{person}{B.~F. Smith}, \bibinfo{person}{P.~E.
  Bjø~rstad}, {and} \bibinfo{person}{W.~D. Gropp}.}
  \bibinfo{year}{1996}\natexlab{}.
\newblock \bibinfo{booktitle}{\emph{Domain decomposition}}.
\newblock \bibinfo{publisher}{Cambridge University Press, Cambridge}.
\newblock
\showISBNx{978-0-521-49589-9}
\urldef\tempurl%
\url{https://mathscinet.ams.org/mathscinet-getitem?mr=1410757}
\showURL{%
\tempurl}


\bibitem[St-Cyr et~al\mbox{.}(2007)]%
        {st-cyr_optimized_2007}
\bibfield{author}{\bibinfo{person}{A. St-Cyr}, \bibinfo{person}{M.~J. Gander},
  {and} \bibinfo{person}{S.~J. Thomas}.} \bibinfo{year}{2007}\natexlab{}.
\newblock \showarticletitle{Optimized multiplicative, additive, and restricted
  additive {Schwarz} preconditioning}.
\newblock \bibinfo{journal}{\emph{SIAM Journal on Scientific Computing}}
  \bibinfo{volume}{29}, \bibinfo{number}{6} (\bibinfo{year}{2007}),
  \bibinfo{pages}{2402--2425}.
\newblock
\showISSN{1064-8275}
\urldef\tempurl%
\url{https://doi.org/10.1137/060652610}
\showDOI{\tempurl}


\bibitem[Toselli and Widlund(2005)]%
        {toselli_domain_2005}
\bibfield{author}{\bibinfo{person}{A. Toselli} {and} \bibinfo{person}{O.
  Widlund}.} \bibinfo{year}{2005}\natexlab{}.
\newblock \bibinfo{booktitle}{\emph{Domain decomposition methods—algorithms
  and theory}}. \bibinfo{series}{Springer {Series} in {Computational}
  {Mathematics}}, Vol.~\bibinfo{volume}{34}.
\newblock \bibinfo{publisher}{Springer-Verlag, Berlin}.
\newblock
\showISBNx{978-3-540-20696-5}
\urldef\tempurl%
\url{https://doi.org/10.1007/b137868}
\showDOI{\tempurl}


\bibitem[Trott et~al\mbox{.}(2022)]%
        {9485033}
\bibfield{author}{\bibinfo{person}{C.~R. Trott}, \bibinfo{person}{D.
  Lebrun-Grandié}, \bibinfo{person}{D. Arndt}, \bibinfo{person}{J. Ciesko},
  \bibinfo{person}{V. Dang}, \bibinfo{person}{N. Ellingwood},
  \bibinfo{person}{R. Gayatri}, \bibinfo{person}{E. Harvey},
  \bibinfo{person}{D.~S. Hollman}, \bibinfo{person}{D. Ibanez},
  \bibinfo{person}{N. Liber}, \bibinfo{person}{J. Madsen}, \bibinfo{person}{J.
  Miles}, \bibinfo{person}{D. Poliakoff}, \bibinfo{person}{A. Powell},
  \bibinfo{person}{S. Rajamanickam}, \bibinfo{person}{M. Simberg},
  \bibinfo{person}{D. Sunderland}, \bibinfo{person}{B. Turcksin}, {and}
  \bibinfo{person}{J. Wilke}.} \bibinfo{year}{2022}\natexlab{}.
\newblock \showarticletitle{Kokkos 3: Programming Model Extensions for the
  Exascale Era}.
\newblock \bibinfo{journal}{\emph{IEEE Transactions on Parallel and Distributed
  Systems}} \bibinfo{volume}{33}, \bibinfo{number}{4} (\bibinfo{year}{2022}),
  \bibinfo{pages}{805--817}.
\newblock
\urldef\tempurl%
\url{https://doi.org/10.1109/TPDS.2021.3097283}
\showDOI{\tempurl}


\bibitem[Wheeler et~al\mbox{.}(2020)]%
        {WheWiLee20}
\bibfield{author}{\bibinfo{person}{M.~F. Wheeler}, \bibinfo{person}{T. Wick},
  {and} \bibinfo{person}{S. Lee}.} \bibinfo{year}{2020}\natexlab{}.
\newblock \showarticletitle{{IPACS: Integrated Phase-Field Advanced Crack
  Propagation Simulator. An adaptive, parallel, physics-based-discretization
  phase-field framework for fracture propagation in porous media}}.
\newblock \bibinfo{journal}{\emph{Computer Methods in Applied Mechanics and
  Engineering}}  \bibinfo{volume}{367} (\bibinfo{year}{2020}),
  \bibinfo{pages}{113124}.
\newblock
\showISSN{0045-7825}
\urldef\tempurl%
\url{https://doi.org/10.1016/j.cma.2020.113124}
\showDOI{\tempurl}


\bibitem[Wick(2011)]%
        {Wi11}
\bibfield{author}{\bibinfo{person}{T. Wick}.} \bibinfo{year}{2011}\natexlab{}.
\newblock \showarticletitle{Fluid-Structure Interactions using Different Mesh
  Motion Techniques}.
\newblock \bibinfo{journal}{\emph{Comput. \& Structures}} \bibinfo{volume}{89},
  \bibinfo{number}{13-14} (\bibinfo{year}{2011}), \bibinfo{pages}{1456--1467}.
\newblock


\bibitem[Wick(2013)]%
        {Wi13_fsi_with_deal}
\bibfield{author}{\bibinfo{person}{T. Wick}.} \bibinfo{year}{2013}\natexlab{}.
\newblock \showarticletitle{Solving Monolithic Fluid-Structure Interaction
  Problems in Arbitrary {L}agrangian {E}ulerian Coordinates with the deal.{II}
  Library}.
\newblock \bibinfo{journal}{\emph{Archive of Numerical Software}}
  \bibinfo{volume}{1} (\bibinfo{year}{2013}), \bibinfo{pages}{1--19}.
\newblock
\urldef\tempurl%
\url{https://doi.org/10.11588/ans.2013.1.10305}
\showDOI{\tempurl}


\bibitem[Wick and Wollner(2021)]%
        {WiWo21}
\bibfield{author}{\bibinfo{person}{T. Wick} {and} \bibinfo{person}{W.
  Wollner}.} \bibinfo{year}{2021}\natexlab{}.
\newblock \showarticletitle{Optimization with nonstationary, nonlinear
  monolithic fluid-structure interaction}.
\newblock \bibinfo{journal}{\emph{Internat. J. Numer. Methods Engrg.}}
  \bibinfo{volume}{122}, \bibinfo{number}{19} (\bibinfo{year}{2021}),
  \bibinfo{pages}{5430--5449}.
\newblock
\urldef\tempurl%
\url{https://doi.org/10.1002/nme.6372}
\showDOI{\tempurl}


\bibitem[Wu and Cai(2011)]%
        {wu_parallel_2011}
\bibfield{author}{\bibinfo{person}{Y. Wu} {and} \bibinfo{person}{X.-C. Cai}.}
  \bibinfo{year}{2011}\natexlab{}.
\newblock \showarticletitle{A parallel two-level method for simulating blood
  flows in branching arteries with the resistive boundary condition}.
\newblock \bibinfo{journal}{\emph{Computers \& Fluids}} \bibinfo{volume}{45},
  \bibinfo{number}{1} (\bibinfo{date}{June} \bibinfo{year}{2011}),
  \bibinfo{pages}{92--102}.
\newblock
\showISSN{0045-7930}
\urldef\tempurl%
\url{https://doi.org/10.1016/j.compfluid.2010.11.015}
\showDOI{\tempurl}


\bibitem[Yamazaki et~al\mbox{.}(2023)]%
        {yamazaki_experimental_2023}
\bibfield{author}{\bibinfo{person}{I. Yamazaki}, \bibinfo{person}{A. Heinlein},
  {and} \bibinfo{person}{S. Rajamanickam}.} \bibinfo{year}{2023}\natexlab{}.
\newblock \showarticletitle{An {Experimental} {Study} of {Two}-level {Schwarz}
  {Domain}-{Decomposition} {Preconditioners} on {GPUs}}.
  \bibinfo{publisher}{IEEE Computer Society}, \bibinfo{pages}{680--689}.
\newblock
\showISBNx{9798350337662}
\urldef\tempurl%
\url{https://doi.org/10.1109/IPDPS54959.2023.00073}
\showDOI{\tempurl}


\bibitem[Zolgharni et~al\mbox{.}(2009)]%
        {Art:Zolgharni:MIT:2009}
\bibfield{author}{\bibinfo{person}{M. Zolgharni}, \bibinfo{person}{P.~D.
  Ledger}, {and} \bibinfo{person}{H.~J. Griffiths}.}
  \bibinfo{year}{2009}\natexlab{}.
\newblock \showarticletitle{Forward Modelling of Magnetic Induction Tomography:
  A Sensitivity Study for Detecting Haemorrhagic Cerebral Stroke}.
\newblock \bibinfo{journal}{\emph{Medical \& Biological Engineering \&
  Computing}} \bibinfo{volume}{47}, \bibinfo{number}{12}
  (\bibinfo{year}{2009}), \bibinfo{pages}{1301--1313}.
\newblock
\showISSN{0140-0118, 1741-0444}
\urldef\tempurl%
\url{https://doi.org/10.1007/s11517-009-0541-1}
\showDOI{\tempurl}


\end{thebibliography}

\end{document}